\documentclass[a4paper]{article}
\usepackage[textwidth=16cm]{geometry}
\usepackage[utf8]{inputenc}
\usepackage{amsmath,amssymb,amsthm}
\usepackage{graphicx}
\usepackage[pdfborder= 0 0 0,citecolor=black,linkcolor=black,colorlinks=True,bookmarksopen=true]{hyperref}
\usepackage{caption}
\usepackage{subcaption}
\usepackage{algorithm}
\usepackage{algorithmicx}
\usepackage{algpseudocode}
\usepackage{tikz}
\usepackage{pgfplots}\pgfplotsset{compat=1.7}
\usepackage{authblk}
\usepackage[square]{natbib}
\usepackage[pagewise]{lineno}
\usepackage{environ}
\usepackage{nag}
\usepackage[colorinlistoftodos,prependcaption,textsize=tiny]{todonotes}
\usepackage{xspace} 
\usepackage{bm}
\usepackage{wasysym}
\usepackage{verbatim}

\usetikzlibrary{plotmarks}
\usepgfplotslibrary{colorbrewer}
\usetikzlibrary{colorbrewer}

\author[1]{P.A. Browne}
\author[2,3]{J. Prettyman}
\author[1]{H. Weller}
\author[2]{T. Pryer}
\author[4]{J. Van lent}

\affil[1]{Department of Meteorology, University of Reading, UK}
\affil[3]{National Physical Laboratory, UK}
\affil[2]{Department of Mathematics, University of Reading, UK}
\affil[4]{FET -- Engineering, Design and Mathematics, University of the West of England, UK}
\affil[*]{Correspondence to \href{mailto:p.browne@reading.ac.uk}{p.browne@reading.ac.uk}}
\title{
Nonlinear solution techniques for solving a Monge-Amp{\`e}re equation for redistribution of a mesh
}

\makeatletter \AtBeginDocument{ \hypersetup{pdftitle= {\@title},pdfauthor= {\@author}}} \makeatother

\date{\today}

\pgfkeys{/pgfplots/legend entry/.code=\addlegendentry{#1}} 

\newcommand{\paper}{article\xspace}
\newcommand{\ma}{Monge-Amp{\`e}re\xspace}
\DeclareMathOperator*{\argmin}{arg\,min}

\makeatletter
\begin{document}
\maketitle
\begin{abstract}
A \ma equation arises when seeking an optimally transported mesh that equidistributes a given monitor function in Cartesian space. This \ma equation is a fully nonlinear PDE, with a source term that is a function of the gradient of the solution. This nonlinear source term is an additional computational challenge that has received little attention from \ma applications in other fields. There are two major components needed to find a solution to the \ma equation: a spatial discretisation and an algorithm to find a solution of the resulting nonlinear algebraic equations. There have been a number of different approaches proposed in the literature to solve the \ma equation but none of which perform consistent comparisons across both algorithmic and discretisation differences.
In this study we explore different algorithmic methods for the \ma equation all within the context of a finite volume spatial discretisation.
We introduce a new linearisation of the \ma equation that neglects the nonlinearities arising from the source term and show that it leads to a method that is fast, robust and free of tuning parameters.
We present numerical experiments that show methods based on this linearisation of the \ma equation are more computationally efficient than those that rely on other techniques such as a parabolic relaxation. Further, the equations resulting from a full linearisation of the \ma equation, equivalent to using Newton's method, can be seen as analogous to an advection-diffusion equation. This allows many tools that exist for computational fluid dynamics to be re-factored easily to solve the \ma equation.
The robustness and efficiency of the newly introduced method gives hope that an adaptive solver for geophysical flows using mesh redistribution can be computationally feasible in the near future. 

\end{abstract}

\section*{Graphical abstract}
\begin{tikzpicture}
\begin{semilogyaxis}[width=10cm,height=8cm,legend style={at={(1.01,0.5)},anchor=west},xlabel=Iteration number,ylabel=Residual of \ma equation,enlargelimits=false]
\addplot [black,mark=.] table[x expr=\coordindex, y index={0}] {_home_pbrowne_OpenFOAM_results_bell_AFP_60_equi};\addlegendentry{Adaptive linearisation};
\addplot [YlOrBr-9-9,mark=.] table[x expr=\coordindex, y index={0}] {_home_pbrowne_OpenFOAM_results_bell_FP2D_60_2.5474_equi};\addlegendentry{Laplacian approx. $\gamma=2.5474$};
\addplot [YlOrBr-9-8,mark=.] table[x expr=\coordindex, y index={0}] {_home_pbrowne_OpenFOAM_results_bell_FP2D_60_2.55_equi};\addlegendentry{Laplacian approx. $\gamma=2.55$};
\addplot [YlOrBr-9-7,mark=.] table[x expr=\coordindex, y index={0}] {_home_pbrowne_OpenFOAM_results_bell_FP2D_60_2.60_equi};\addlegendentry{Laplacian approx. $\gamma=2.60$};
\addplot [YlOrBr-9-6,mark=.] table[x expr=\coordindex, y index={0}] {_home_pbrowne_OpenFOAM_results_bell_FP2D_60_2.65_equi};\addlegendentry{Laplacian approx. $\gamma=2.65$};
\addplot [YlOrBr-9-5,mark=.] table[x expr=\coordindex, y index={0}] {_home_pbrowne_OpenFOAM_results_bell_FP2D_60_2.70_equi};\addlegendentry{Laplacian approx. $\gamma=2.70$};
\addplot [YlOrBr-9-4,mark=.] table[x expr=\coordindex, y index={0}] {_home_pbrowne_OpenFOAM_results_bell_FP2D_60_2.75_equi};\addlegendentry{Laplacian approx. $\gamma=2.75$};
\addplot [BuPu-9-3,mark=.] table[x expr=\coordindex, y index={0}] {_home_pbrowne_OpenFOAM_results_bell_FP2D_60_2.80_equi};\addlegendentry{Laplacian approx. $\gamma=2.80$};
\addplot [BuPu-9-4,mark=.] table[x expr=\coordindex, y index={0}] {_home_pbrowne_OpenFOAM_results_bell_FP2D_60_2.85_equi};\addlegendentry{Laplacian approx. $\gamma=2.85$};
\addplot [BuPu-9-5,mark=.] table[x expr=\coordindex, y index={0}] {_home_pbrowne_OpenFOAM_results_bell_FP2D_60_2.90_equi};\addlegendentry{Laplacian approx. $\gamma=2.90$};
\addplot [BuPu-9-6,mark=.] table[x expr=\coordindex, y index={0}] {_home_pbrowne_OpenFOAM_results_bell_FP2D_60_2.95_equi};\addlegendentry{Laplacian approx. $\gamma=2.95$};
\addplot [BuPu-9-7,mark=.] table[x expr=\coordindex, y index={0}] {_home_pbrowne_OpenFOAM_results_bell_FP2D_60_3.00_equi};\addlegendentry{Laplacian approx. $\gamma=3.00$};
\addplot [BuPu-9-8,mark=.] table[x expr=\coordindex, y index={0}] {_home_pbrowne_OpenFOAM_results_bell_FP2D_60_3.05_equi};\addlegendentry{Laplacian approx. $\gamma=3.05$};
\addplot [BuPu-9-9,mark=.] table[x expr=\coordindex, y index={0}] {_home_pbrowne_OpenFOAM_results_bell_FP2D_60_3.10_equi};\addlegendentry{Laplacian approx. $\gamma=3.10$};

\addplot [red,mark=.,solid] table[x expr=\coordindex, y index={0}] {_home_pbrowne_OpenFOAM_results_bell_NEWTON2D-vectorGradc_m_60_equi};\addlegendentry{Full linearisation (Newton's method)};
\end{semilogyaxis}
\end{tikzpicture}

\noindent \textbf{Keywords:} \ma equation, mesh redisribution, r-adaptivity, fixed point iterations, optimal transport, Newton's method

\clearpage{}\section{Introduction}\label{sec:intro}
Variable resolution meshes can be advantageous for the numerical solution of PDEs when variations or sensitivity to errors are greater in some areas than others. This occurs, for example, for the numerical prediction of tropical cyclones or for regional weather forecasts \citep{Wang2001,Piani2000}.
It may be advantageous that the mesh vary through time
tracking 
atmospheric fronts \citep{Kuhnlein2012} or tsunamis \citep{Harig2008}.
Adding (or subtracting) mesh points where more (or fewer) points are required 
is known as \emph{h-adaptivity}. With this type of adaptivity, the connectivity of the mesh (and possibly the total number of points defining the mesh) can change through time \citep{Kimura2013}.

R-adaptivity (mesh redistribution) involves keeping the mesh connectivity fixed but moving the mesh points. R-adaptivity retains fixed data structures associated with the mesh, need not entail mapping the solution from one mesh to another, may not lead to load balancing problems on parallel computers, can be easily incorporated into legacy code and can be designed to give meshes that vary smoothly in space and time. R-adaptivity can lead to smoothly graded meshes \citep{Budd2009,Cao2003}, which are may help to alleviate wave reflections or other errors associated with abrupt changes in resolution \citep{Guba2014,Long2011,Vichnevetsky1987}.

\emph{Optimal transport} is a good technique for mesh-redistribution because it guarantees to find a mesh which is equidistributed with respect to a monitor function that is not tangled \citep{Budd2009a,Weller2016}.
We aim to find a technique for finding optimally transported meshes which is robust, free of tunable parameters and fast enough to be used in numerical weather prediction when the mesh is being moved every time-step.
In this \paper we will compare some of the existing solution techniques and introduce a new, faster and robust method.

A number of spatial discretistations and algorithms for solving non-linear equations have been introduced and used for solving the \ma equation \citep[e.g.][]{Budd2009a,Chacon2011,Weller2016} but some of these have free parameters and it is not clear which of these might be fast and robust enough for frequent adaptations during numerical weather prediction. 
In this \paper we shall use the finite volume spatial discretisation and perform a direct comparison of a number of existing and new algorithms for solving the \ma equation in the context of mesh redistribution.
\autoref{sec:background} introduces the mathematical background of using the \ma equation to find an optimally transported mesh as well as summarising the literature on solving the \ma equation. In \autoref{sec:tech} we describe existing algorithms for solving the \ma equation, and introduce new algorithms based on linearisations of different terms of the equation. In \autoref{sec:num} we describe numerical experiments and diagnostics that we use to test the various algorithms. Finally, in \autoref{sec:conc} we draw conclusions about the different algorithmic methods for solving the \ma equation for mesh redistribution.

\section{Mathematical background}\label{sec:background}
This section describes the theory of optimal transport for mesh redistribution and surveys the various approaches that have been taken in the literature to solve the resulting nonlinear problem.

\begin{figure*}[h]
\begin{subfigure}{0.48\textwidth}
\includegraphics[width=\textwidth]{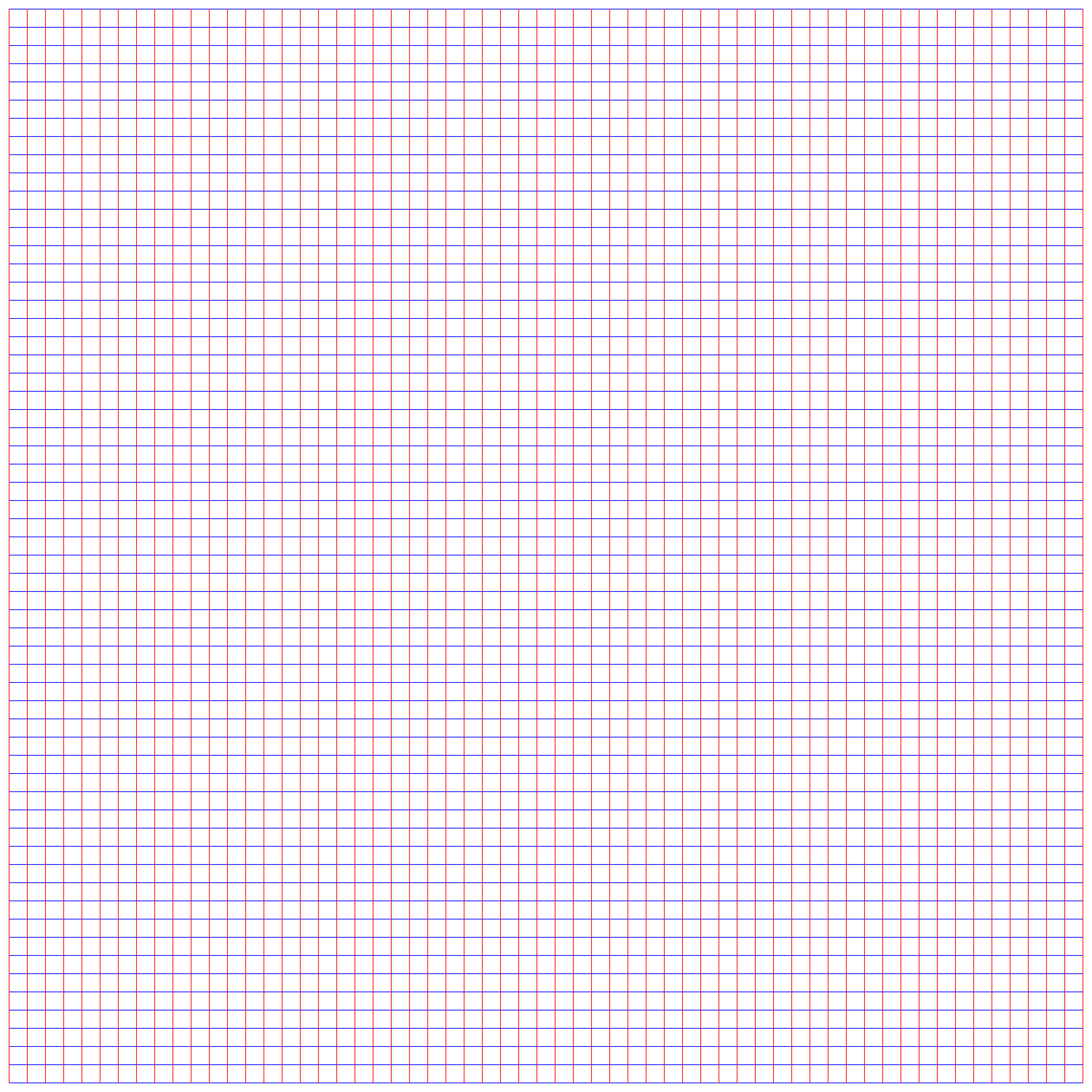}
\caption{Original computational mesh $\mathcal{T}_c$}\label{fig:example_comp}
\end{subfigure}
\hfill
\begin{subfigure}{0.48\textwidth}
\includegraphics[width=\textwidth]{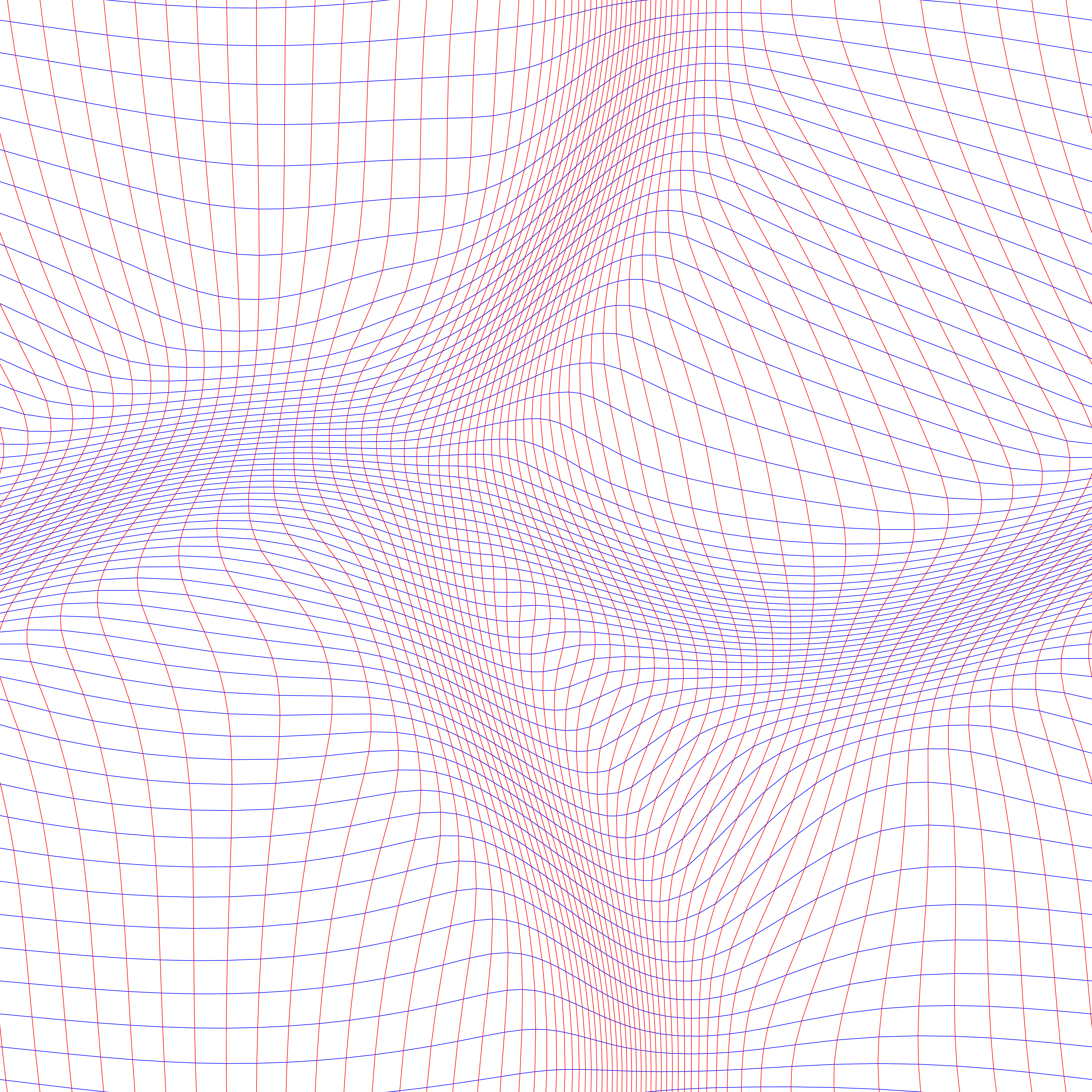}
\caption{Adapted physical mesh $\mathcal{T}_p$}
\end{subfigure}
\caption{The action of a map $\bm{f}$ takes $\mathcal{T}_c$ to $\mathcal{T}_p$}\label{fig:example}
\end{figure*}

Moving meshes, or r-adaptivity, is concerned with relocating the points of a mesh while keeping the mesh connectivity fixed. That is, an original mesh is given, $\mathcal{T}_c$, and it is transformed under a map $\bm{f}$ to $\mathcal{T}_p$ such that $\mathcal{T}_p = \bm{f}(\mathcal{T}_c)$. \autoref{fig:example} shows an example of these two meshes.
The map, $\bm{f}$, needs to move points to regions which require higher resolution whilst retaining
certain desirable properties of the original mesh, such as skewness, orthogonality or convexity.

\subsection{Monitor function equidistribution}
To help with notation, we introduce the computational and physical spaces, $\Omega_c$ and $\Omega_p$, respectively. We have $\mathcal{T}_c \in \Omega_c$ and $\mathcal{T}_p \in \Omega_p$ and we write $\bm{\xi} \in \mathcal{T}_c$ and $\bm{x} \in \mathcal{T}_p$ to distinguish between elements of the two spaces.

To control the mesh adaptivity we assume that we have a scalar valued \emph{monitor function} $m(\bm{x})>0$ given \textit{a priori} which is large in parts of the domain which require a dense mesh and small in the parts of the domain which can be rarefied.
We shall try to \emph{equidistribute} this monitor function under the desired map. That is, for each arbitrary volume in the physical mesh, the total monitor function contained within that volume is proportional to the volume of the inverse image of the volume, i.e. the  corresponding volume in the original computational mesh. 
This can be rigorously written \citep{Budd2009,huang2011,Weller2016} as the equidistribution equation
\begin{equation}\label{eqn:equi}
m(\bm{x})|J(\bm{\xi})| = c,
\end{equation}
where $J$ is the Jacobian of the map $\bm{f}$ and $c$ is a constant. If $\Omega_c$ and $\Omega_p$ are one-dimensional intervals then the solution of \eqref{eqn:equi} uniquely defines the map $\bm{f}$. In higher dimensional problems (such as we consider in this \paper) the equidistribution equation is not enough to uniquely define $\bm{f}$; for example in a problem involving symmetries, $\bm{f}$ could be rotated and still satisfy equidistribution.

\subsection{Optimally transported meshes}
In order to arrive at a well-posed problem to find the adapted mesh, we must choose a method of regularising the equidistribution equation. There are many possible ways of doing so, such as specifying local conditions for the mesh. However, we note that the original computational mesh is likely to have desirable qualities, such as orthogonality, regularity etc. 
Hence it appears reasonable to seek an equidistributed mesh that is as close as possible to the computational mesh. In this sense we want to find an optimal map $\bm{f}^*$ such that
\begin{equation}
\bm{f}^*=\argmin_{\bm{f}: \Omega_c \to \Omega_p} \int_{\Omega_c} |\bm{\xi} - \bm{f}(\bm{\xi})|^2\ \mathrm{d}\bm{\xi}.
\end{equation}
Once this map has been found, we call $\bm{f}^*$ an \emph{optimal transport map} and the resulting mesh $\mathcal{T}_p = \bm{f}^*(\mathcal{T}_c)$ an \emph{optimally transported mesh}.

\subsection{The \ma equation}
A remarkable result from optimal transport theory is that the optimal mapping $\bm{f}^*$ is unique and can be written as the gradient of a convex potential \citep{Brenier1991}. Hence
\begin{equation}\label{eqn:pot}
\bm{x} = \bm{f}^*(\bm{\xi}) = \nabla P(\bm{\xi})
\end{equation}
where $P$ is a convex mesh potential. As $P$ is convex it is easy to show that an optimally transported mesh does not exhibit tangling \citep{Budd2009a}. Substituting \eqref{eqn:pot} into \eqref{eqn:equi} we arrive at a Monge-Amp{\`e}re equation:
\begin{equation}\label{eqn:ma1}
m(\nabla P)|H(P)| = c
\end{equation}
where $H(P)$ is the Hessian of the potential $P$, or equivalently the Jacobian of the map $\bm{f}^*$. Writing
\begin{equation}
\phi(\bm{\xi}) = P - \frac{|\bm{\xi}|^2}{2}
\end{equation}
we arrive at the equation
\begin{equation}\label{eqn:ma2}
m(\bm{\xi} + \nabla \phi)|I+H(\phi)| = c
\end{equation}
which is a fully nonlinear elliptic PDE in the new mesh potential $\phi$. The desired physical mesh is therefore given as
\begin{equation}
\bm{x} = \bm{\xi} + \nabla \phi.
\end{equation}
As $m>0$, we can re-write the \ma equation \eqref{eqn:ma2} as
\begin{equation}\label{eqn:ma3}
| I + H(\phi)| = \frac{c}{m(\bm{x})} .
\end{equation}

To numerically solve the \ma equation requires two components: a spatial discretisation to get a discrete representation of the equation and an algorithm which finds the solution to the resulting system of nonlinear algebraic equations.

\subsubsection{Spatial Discretisations\label{sub:introSpace}}

Spatial discretisation of the \ma equation has mostly focused on
the discretisation of $|H(\phi)|$ \citep[e.g.][]{Froese2011}. \citet{Oberman2008}
developed wide stencil finite difference techniques 
based on the \citet{Barles1991} framework 
to ensure discrete monotonicity guaranteeing convergence of the
numerical solution to the unique convex viscosity solution. 
Mixed
finite-element methods have been used by \citet{Dean2006,Feng2009,Lakkis2013,Neilan2014}. Finite differences have been used extensively in 
solving the \ma equation for mesh redistribution problems \citep[e.g.][]{Delzanno2008,Sulman2011}.
\citet{Budd2009,Browne2014a}
used simple finite differences with filtering of the right hand side
and smoothing of the Hessian in order to ensure that a convex solution
is found. \citet{Weller2016} used finite volumes to solve the \ma equation
on a plane and on the sphere and explored a number of techniques for
solving the $m(\bm{\xi}+\nabla\phi)$ term on the right hand side
of eqn \eqref{eqn:ma3}
while \citet{Froese2012} used wide stencil finite differences to discretise
$\nabla\phi$ on the right hand side of the \ma equation and \citet{Saumier2015} experimented with finite differences and a spectral method.
\citet{Feng2016} give a reformulation of the \ma problem as a Hamilton-Jacobi-Bellman equation which removes the constraint of convexity on the solution as a by-product of this reformulation.

In this \paper we shall consider using only the finite volume technique similarly to \citet{Weller2016} in order to compare various differences in algorithmic processes for solving the \ma equation.
 
\subsubsection{Non-linear Equation Solution Algorithms\label{sub:introAlgorithm}}

The \ma equation can be posed as the nonlinear PDE
\[|H(\phi)| = g.\]

The majority of the literature on solving this \ma equation is devoted to the case where the source term $g$ is a given function of spatial location on the computational mesh \citep[e.g.][]{Dean2006,Feng2009,Oberman2008}. For the mesh redistribution, this source term becomes a function of the solution $\phi$ itself (i.e. $g= g(\nabla \phi)$) since $m$ is defined in physical space. Techniques for treating this nonlinearity in the source term have been investigated by \citet{Froese2012} and \citet{Weller2016}. This nonlinearity also occurs in the case of prescribed Gauss curvature \citep[see for example][]{Pryer2012}.

\citet{Benamou2010} give a concise overview of the algorithmic methods for solving the \ma equation when $g$ is not directly a function of the solution. \citet{Dean2006a,Dean2006} use an Augmented Lagrangian method to the \ma equation after it has been reformulated as a saddle point problem. \citet{Dean2006} use a nonlinear least-squares solution to minimise an appropriate functional of the solution. \citet{Benamou2010} use a fixed point method derived by taking a Poisson approximation to the \ma equation. A related fixed point method using a Poisson approximation was introduced by \citet{Weller2016} for solving an equation of \ma type for mesh redistribution on a spherical manifold.

For mesh redistribution via optimal transport in Cartesian geometry, \citet{Delzanno2008} and \citet{Chacon2011} use a multigrid-preconditioned inexact Newton-Krylov method with damping. A parabolic relaxation to the \ma equation was proposed by \citet{Budd2009a} and has been used successfully in a number of studies (\citep{Budd2013,Browne2014a}), with other forms of parabolic relaxation investigated by \citet{Sulman2011}.

\clearpage{}

\section{Details of techniques for solving the \ma equation}\label{sec:tech}

In  this  section  we  shall  describe in detail algorithms  for  solving  the \ma equation  in  the  context  of mesh redistribution where the source term is a function of the gradient of the solution.  In particular, we will introduce a linearisation about the current approximate solution which will lead to a fixed point method with no free parameters.
We will then describe  the  spatial  discretisation.The essential computational components of the algorithms are summarised in Table \ref{table:comp}.

\subsection{Fixed point iterations}\label{sec:fixed_point}
In order to find a solution to the \ma equation, \citet{Weller2016} introduced a fixed point method, which we describe here.

Firstly, as $m>0$, we can re-write the \ma equation \eqref{eqn:ma2} as
\begin{equation}
| I + H(\phi)| - \frac{c}{m(\bm{x})} = 0
\end{equation}
Let $\phi^n$ denote the approximation to the solution at iteration $n$. Taking a Taylor's series expansion of $|I+H(\phi^{n+1})|$ about $\phi^0=0$ we see
\begin{equation}\label{alabama}
|I+H(\phi^{n+1})| = 1 + \nabla^2 \phi^{n+1} + \mathcal{N}(\phi^{n+1})
\end{equation}
where $\mathcal{N}(\phi^{n+1})$ are higher order, nonlinear terms. In 2D we have that the nonlinear terms are precisely $\mathcal{N}(\phi^{n+1}) = |H(\phi^{n+1})|$.

We want to solve
\begin{equation}\label{colorado}
|I+H(\phi^{n+1})| = \frac{c}{m(\bm{x}^{n+1})}.
\end{equation}
Substituting \eqref{alabama} into \eqref{colorado}, and subtracting $|I+H(\phi^n)|$ we obtain
\begin{equation}
1+\nabla^2 \phi^{n+1} +\mathcal{N}(\phi^{n+1})- \mathcal{N}(\phi^n) = 1+\nabla^2 \phi^{n} - |I+H(\phi^{n})|+\frac{c}{m(\bm{x}^{n+1})}.
\end{equation}
Dropping the change in nonlinear terms between iterations $n$ and $n+1$ gives a
fixed point iteration in $\phi^n$ such that, given $\bm{x}^0$ and $\phi^0$,
\begin{equation}\label{eqn:fp1}
1+\nabla^2 \phi^{n+1} = 1 + \nabla^2 \phi^{n} - | I + H(\phi^n)| + \frac{c}{m(\bm{x}^n)}, \qquad \forall n \in \mathbb{N}.
\end{equation}
It is easy to see that a fixed point of this iteration will solve the fully nonlinear \ma equation. For this fixed point iteration to be stable we require that the nonlinear terms we have omitted, $\mathcal{N}(\phi^{n+1})- \mathcal{N}(\phi^n)$, are small. Equivalently, $|H(\phi^{n+1})|\approx|H(\phi^n)|$. Note that these are nonlinear functions of variables the same order as $\phi$.

As \citet{Weller2016} showed (and can be seen in Section \ref{sec:num}) this fixed point iteration is not always convergent. Under-relaxation can force this equation to converge. That is, introducing a scalar constant $\gamma$, we can weight the linear terms in \eqref{eqn:fp1} such that
\begin{equation}\label{eqn:fp11}
\gamma\nabla^2 \phi^{n+1} = \gamma\nabla^2 \phi^{n} - | I + H(\phi^n)| + \frac{c}{m(\bm{x}^n)}, \qquad \forall n \in \mathbb{N}.
\end{equation}
For large $\gamma$ there is less weight given to the error in the nonlinear \ma equation and hence only a small update is made to the approximate solution $\phi$. Vice versa, the smaller $\gamma$, the more the nonlinear problem is treated directly, and hence the large the update from each fixed point iteration. We have no theoretical basis on which to find the optimal value of $\gamma$ which, as will be seen in Section \ref{sec:num}, is dependent on the monitor function in question. In order to have a method with no free parameters, we introduce a different linearisation of the \ma equation.

\subsection{A fixed point method with an adaptive linearisation}\label{sec:lin}
Instead of linearising $|I+H(\phi^{n+1})|$ about $\phi^0$, we can instead linearise about a current iterate $\phi^n$. Writing $\phi^{n+1} = \phi^n+\varepsilon\psi$ it can be shown that
\begin{equation}\label{delaware}
|I+H(\phi^{n+1})| = |I+H(\phi^n)| + \nabla \cdot A^n\nabla \varepsilon \psi + \mathcal{N}(\varepsilon \psi),
\end{equation}
where $A^n$ is the matrix of cofactors of $I+H(\phi)$ and $\mathcal{N}$ is some nonlinear function. In 2D
\begin{equation}
A^n = \begin{bmatrix} 1 + \phi^n_{yy} & -\phi^n_{xy} \\ -\phi^n_{xy} & 1+\phi^n_{xx}\end{bmatrix}
\end{equation}
and $\mathcal{N}(\varepsilon \psi) =  \varepsilon^2 |H(\phi)|$.
In 3D, a more involved computation can show $\mathcal{N}(\varepsilon \psi) =  \varepsilon^3 \mathcal{\tilde{N}}(\psi)$ and 
\begin{equation}
A^n = \begin{bmatrix} 
1 + \phi^n_{yy} + \phi^n_{zz} + \phi^n_{yy}\phi^n_{zz} - \phi^n_{yz}\phi^n_{yz} & -\phi^n_{xy} -\phi^n_{xy}\phi^n_{zz} + \phi^n_{xz}\phi^n_{yz} & -\phi^n_{xz} -\phi^n_{xz}\phi^n_{yy}+\phi^n_{xy}\phi^n_{yz}\\
-\phi^n_{xy} -\phi^n_{xy}\phi^n_{zz} +\phi^n_{xz}\phi^n_{yz} & 1+\phi^n_{xx} + \phi^n_{zz} + \phi^n_{xx}\phi^n_{zz} - \phi^n_{xz}\phi^n_{xz}&-\phi^n_{yz} - \phi^n_{xx}\phi^n_{yz}+\phi^n_{xy}\phi^n_{xz}\\
-\phi^n_{xz} - \phi^n_{xz}\phi^n_{yy} + \phi^n_{xy}\phi^n_{yz} & -\phi^n_{yz} - \phi^n_{xx}\phi^n_{yz} + \phi^n_{xy}\phi^n_{xz} & 1+ \phi^n_{xx}+\phi^n_{yy}+\phi^n_{xx}\phi^n_{yy}-\phi^n_{xy}\phi^n_{xy}
\end{bmatrix}.
\end{equation}

Hence defining $A^n$ as above
and substituting \eqref{delaware} into the \ma equation \eqref{colorado}  we obtain
\begin{equation}\label{florida}
|I+H(\phi^{n})| + \nabla \cdot \left(A^n\nabla \varepsilon \psi\right) + \mathcal{N}(\varepsilon \psi) = \frac{c}{m(\bm{x}^{n+1})}.
\end{equation}
Dropping the 
terms proportional to $\varepsilon^d$ where $d$ is the dimension of the space
gives a
fixed point iteration in $\phi^n$ 
such that, given $\bm{x}^0$ and $\phi^0$,
\begin{equation}\label{eqn:fp2}
\nabla \cdot \left(A^n\nabla \varepsilon \psi\right) = - | I + H(\phi^n)| + \frac{c}{m(\bm{x}^n)}, \qquad \forall n \in \mathbb{N}.
\end{equation}
As in \eqref{eqn:fp1}, a fixed point of \eqref{eqn:fp2} solves the \ma equation. The nonlinear term we have omitted in \eqref{eqn:fp1} is $\mathcal{N}(\varepsilon \psi) = \varepsilon^d|H(\psi)|$. When $\varepsilon$ is small (i.e. $\phi^{n+1}\approx\phi^n$) the nonlinear terms will be smaller than those omitted in the fixed point method \eqref{eqn:fp1}. 
Note that at $\phi^0=0$, the \eqref{eqn:fp2} reduces to the initial fixed point method as $H(0)=0 \implies A^0=0$.

The iterative method for solving the \ma equation given in \eqref{eqn:fp2} can have numerical difficulties when the discretised, matrix equation becomes indefinite. In this case, the ellipticity property of the original \ma equation is lost and the 
nonconvex solutions to the \ma equation can be generated which lead to tangled meshes and numerical divergence.
The indefiniteness of \eqref{eqn:fp2} is caused by the matrix $A^n$ being numerically indefinite. Therefore, we can modify \eqref{eqn:fp2} to maintain ellipticity by

\begin{equation}\label{eqn:fp3}
\nabla \cdot \left(B^n\nabla \varepsilon \psi \right) = - | I + H(\phi^n)| + \frac{c}{m(\bm{x}^n)}, \qquad \forall n \in \mathbb{N},
\end{equation}
where
\begin{equation}\label{eqn:fp3gamma}
B^n = A^n + \gamma I
\end{equation}
and $\gamma$ is defined as
\begin{equation}\label{eqn:eigenvalue}
\gamma := \begin{cases}
0 \qquad &\text{if} \qquad \min\sigma[A^n] > 0\\
\epsilon - \min \sigma[A^n] \qquad &\text{if} \qquad \min\sigma[A^n] \le 0.
\end{cases}
\end{equation}
The constant $\epsilon>0$ is chosen to avoid round-off errors (we have taken $\epsilon=10^{-5}$), and $\sigma[A^n]$ refers to the spectrum of $A^n$. This process simply shifts the eigenvalues of the matrix $A^n$ so that they remain positive. One can check that this choice of $\gamma$ ensures that $B^n \succ 0$. We shall refer to the algorithm given in \eqref{eqn:fp3} as the Adaptive Fixed Point method (AFP).

There are two key remarks about about \eqref{eqn:fp3}: firstly there are no free parameters in the method and, secondly, the 
scalar $\gamma$ 
can be spatially varying. The lack of free parameters makes this very attractive as a robust method for adapting a mesh. Such robust methods will be required for the future coupling of optimally transported meshes into operational forecasting models where the monitor function that can vary in time as well as space is not known \textit{a priori}. 
Instead of using one global $\gamma$ (as for the FP method in Section \ref{sec:fixed_point}) to regularise the global system, this choice of regularisation shifts only the eigenvalues of the local matrices $A^n$ in the parts of the domain where they are numerically indefinite.

A further remark about \eqref{eqn:fp3} is that is can be discretised in space using a finite volume, finite element or finite difference method and solved implicitly for $\phi^{n+1}$. It can be seen as a nonconstant tensor coefficient Poisson equation. This type of equation is particularly difficult to solve with a spectral method as it is nonseparable \citep{boyd2013}.

\subsection{Newton's method for the \ma equation}\label{sec:newton}
Starting from \eqref{colorado}, we wish to solve 
\begin{equation}\label{colorado2}
|I+H(\phi^{n+1})| = \frac{c}{m(\bm{x}^{n+1})}.
\end{equation}

Numerical solution using Newton's method 
involves linearising not only $|I+H(\phi)|$, but also the nonlinear right hand side $\frac{c}{m}(\phi)$.

Linearising the right hand side $\frac{c}{m(\bm{x}^{n+1})}$ about $\bm{x}^n$ gives
\begin{align}
\frac{c}{m}(\bm{x}^{n+1}) &= 
\frac{c}{m}(\bm{x}^{n}) + \nabla_{\bm{x}}\left( \frac{c}{m}(\bm{x}^n)\right) \cdot (\bm{x}^{n+1}-\bm{x}^n).
\end{align}
Recall $\bm{x}^{n+1} = \bm{\xi} + \nabla \phi^{n+1}$ and $\bm{x}^{n} = \bm{\xi} + \nabla \phi^{n}$, thus $\bm{x}^{n+1}-\bm{x}^{n} = \nabla \phi^{n+1} - \nabla \phi^{n}$. Also we write $\phi^{n+1} = \phi^n + \varepsilon \psi$. Hence 
\begin{align}
\frac{c}{m}(\bm{x}^{n+1}) &= 
\frac{c}{m}(\bm{x}^{n}) + \nabla_{\bm{x}}\left( \frac{c}{m}(\bm{x}^n)\right) \cdot (\nabla \phi^{n+1}-\nabla \phi^n)\\
&= \frac{c}{m}(\bm{x}^{n}) + \nabla_{\bm{x}}\left( \frac{c}{m}(\bm{x}^n)\right) \cdot \nabla (\varepsilon \psi).\label{arkansas}\end{align}
This can be incorporated into \eqref{florida} to obtain the linearisation of the full \ma equation. Note that \eqref{arkansas} contains gradients in two different spaces: $\nabla$, the gradient on the computational mesh, and $\nabla_{\bm{x}}$, the gradient 
on the physical mesh.

Substituting \eqref{arkansas} and \eqref{delaware} into \eqref{colorado2} gives:
\begin{equation}
|I+H(\phi^{n})| + \nabla \cdot \left(B^n\nabla \varepsilon\psi\right) = \frac{c}{m(\bm{x}^{n})} + \nabla_{\bm{x}}\left( \frac{c}{m}(\bm{x}^n)\right) \cdot \nabla (\varepsilon\psi).
\end{equation}
Rearranging the terms leads to
\begin{equation}\label{eqn:newton1}
\nabla \cdot \left(B^n\nabla \varepsilon\psi\right) - \nabla_{\bm{x}}\left( \frac{c}{m}(\bm{x}^n)\right) \cdot \nabla (\varepsilon\psi) + |I+H(\phi^{n})| - \frac{c}{m(\bm{x}^{n})} = 0.
\end{equation}

Iterating equation \eqref{eqn:newton1} is Newton's method for solving the \ma equation. Written in this form it can be seen to have a physical interpretation: it is an advection--diffusion equation for $\varepsilon\psi$ with tensorial diffusion coefficient $B^n$, advection velocity $\nabla_{\bm{x}}\left( \frac{c}{m}(\bm{x}^n)\right)$ and source term $-|I+H(\phi^{n})| + \frac{c}{m(\bm{x}^{n})}$.

As we can write Newton's method for solving the \ma equation as an advection--diffusion equation, we can use the tools and knowledge from
computational fluid dynamics to find a numerical solution.
This will be discussed further in Section \ref{newton_numerics}.

\subsection{The Parabolic Monge-Amp{\`e}re method}\label{sec:pma}
The parabolic \ma (PMA) method for solving the \ma equation was introduced by \citet{Budd2009a}. They consider a pseudo-time equation
\begin{equation}\label{eqn:pma1}
(I-\gamma \nabla^2)\phi_\tau = \left[m(\bm{x})|I+H(\phi)|\right]^{\tfrac{1}{d}}
\end{equation}
where $\gamma$ is a smoothing coefficient (constant over space), $\phi_\tau$ is the pseudo-time derivative of the potential $\phi$ and $d$ is the dimension of the space in which the \ma equation is being solved. \citet{Budd2009a} showed that as $\tau \to \infty$, $\nabla \phi \to \nabla P$ where $\nabla P$ solves the original \ma equation. Other parabolic relaxations of the \ma equation have been investigated in the literature. For example, \citet{Sulman2011} use a logarithmic form of relaxation instead of the power law in \eqref{eqn:pma1}.

In the PMA method, only the gradient of the potential is of interest. From a numerical perspective, it is useful to note that the potential is always increasing: $m>0$ and $|H|>0$ as $H$ is symmetric positive definite. Hence this potential can be modified by the addition of any constant without affecting the PMA method's convergence to the solution of the original \ma equation.

We can choose this constant to be
\begin{equation}
c = -\int_{\Omega} (I-\gamma \nabla^2)^{-1}\left[m(\bm{x})|I+H(\phi)|\right]^{\tfrac{1}{d}}
\end{equation}
so that
\begin{equation}
\int_{\Omega} \phi_\tau = 0 \implies \int_{\Omega} \phi^n = \int_{\Omega} \phi^0 \quad \forall n \in \mathbb{N}.
\end{equation}

With this choice of constant, PMA can be written as a fixed point method in $\phi^n$
\begin{equation}
(I-\gamma \nabla^2)\phi^{n+1} = (I-\gamma \nabla^2)\phi^n + \delta t \left[m(\bm{x}^n)|I+H(\phi^n)|\right]^{\tfrac{1}{d}} + c.
\end{equation}

Note that PMA has 2 separate parameters to choose in order to solve the \ma equation: $\gamma$ and $\delta t$. We consider the PMA method here as a benchmark with which to compare the fixed point and adaptive linearisation methods due to the considerable literature using PMA for mesh adaptation (e.g. \citep{Budd2009a,Budd2013,Browne2014a} etc).

The algorithmic methods which we consider in this \paper are summarised in Table \ref{table:comp} that shows what matrix equation is solved at each timestep and what free parameters each method needs to set.
 
\begin{table}[]
\centering
\caption{Comparison of algorithmic methods for solving the \ma equation}
\label{table:comp}
\begin{tabular}{|l|l|l|l|}
\hline
Method                     & Inversion                 & Free parameters      & Advection term \\ \hline
Fixed point (FP)           & $\nabla^2$                & $\gamma$             & None           \\
Adaptive fixed point (AFP) & $\nabla \cdot B^n \nabla$ & None                 & None           \\
Newton's method            & $\nabla \cdot B^n \nabla$ & $\delta$             & Yes            \\
Parabolic relaxation (PMA) & $I-\gamma \nabla^2$       & $\gamma$, $\delta t$ & None           \\ \hline
\end{tabular}
\end{table}

\subsection{Spatial Discretisation\label{sub:space}}

We use a finite volume technique to discretise the \ma equation in
space using OpenFOAM \citep{openfoam}, following \citet{Weller2016}.
The discretisation assumes that all finite volume cells are three dimensional with two dimensional faces between cells. The test cases in this paper use one layer of cells to represent the two dimensional domain.
The prognostic variable $\phi$ is stored at the centre of each cell. To calculate the value of the Laplacian in cell $i$, $(\nabla^2 \phi)_i$ we use the Divergence Theorem, and write it as a sum over each face $f$ of the cell. Hence
\begin{equation}\label{eqn:lap1}
(\nabla^2 \phi)_i \approx \frac{1}{V_i} \sum_{f\in i} |\bm{S}_f| \frac{\phi_{i_f} - \phi_i}{|\bm{d}_f|}
\end{equation}
where $V_i$ is the volume of cell $i$, $\bm{S}_f$ is the vector normal to the face $f$ with magnitude $|\bm{S}_f|$ equal to the area of the face $f$, $i_f$ refers to the cell connected to cell $i$ via the face $f$ and $|\bm{d}_f|$ is the distance between the centre of cell $i_f$ and the centre of cell $i$. The notation $f\in i$ refers to a face $f$ of cell $i$.

The gradient of $\phi$ at cell centres. $\nabla_c \phi$ is computed using the divergence theorem such that
\begin{equation}\label{eqn:div_thm}
(\nabla_c \phi)_i = \frac{1}{V_i}\sum_{f\in i} \phi_f \bm{S}_f
\end{equation}
where $\phi_f$ is the value of $\phi$ linearly interpolated onto face $f$.
The gradient of $\phi$ can be computed in the direction normal to a face $f$ by 
\begin{equation}\label{maine} \nabla_{nf} \phi = \frac{\phi_{i_f} - \phi_i}{|\bm{d}_f|}. \end{equation}

To calculate the value of the Hessian in cell $i$, $H_i(\phi)$ we first compute $\nabla_c \phi$ which is valid on the cell centres and linearly interpolate this onto the faces to give $\widetilde{\nabla_f \phi}$. The normal component of the full gradient on the face, $\nabla_f \phi$, is corrected such that
\begin{equation}
\nabla_f \phi = \widetilde{\nabla_f \phi} + \left(\nabla_{nf} \phi - (\widetilde{\nabla_f \phi} \cdot \hat{\bm{S}}_f)\right) \hat{\bm{S}}_f
\end{equation}
where $\hat{\bm{S}}_f$ is the unit vector normal to face $f$.

Then we compute the gradient of this which is valid at the cell centres using the divergence theorem. 
\begin{equation}
H_i(\phi) = \nabla \nabla \phi = \frac{1}{V_i}\sum_{f\in i} (\nabla_f \phi) \bm{S}_f
\end{equation}
We store the mesh points at the cell corners. Therefore to update the mesh we need to have a gradient $\nabla \phi$ valid on cell corners. To do this, we reconstruct the gradient at the corners from $\nabla_{nf} \phi$ which, in the case of a uniform grid, is the average the gradient on each face that contains the cell corner in question.

\subsubsection{Spatial discretisation of the advection term for Newton's method}\label{newton_numerics}
Recall that in \eqref{eqn:newton1} we wrote Newton's method for solving the \ma equation as an advection--diffusion equation. 
This will be solved using existing functions in OpenFOAM which can solve advection-diffusion equations implicitly if the advection term is written in conservative form.
To go from advective to conservative form, we note the following vector calculus identity:
\begin{equation}\label{eqn:identity}
\nabla_{\bm{x}} a \cdot \nabla_{\bm{y}} b = \nabla_{\bm{y}} \cdot ((\nabla_{\bm{x}} a)b) - b \nabla_{\bm{y}} \cdot \nabla_{\bm{x}} a
\end{equation}

Using the identity \eqref{eqn:identity}, the advection term in \eqref{eqn:newton1} becomes
\[
\nabla_{\bm{x}}\left( \frac{c}{m}(\bm{x}^n)\right) \cdot \nabla \psi = \nabla \cdot (\nabla_{\bm{x}}\left( \frac{c}{m}(\bm{x}^n)\right)\varepsilon\psi) - \varepsilon\psi \nabla \cdot \nabla_{\bm{x}} \left( \frac{c}{m}(\bm{x}^n)\right)
\]
and so (changing sign to follow the usual convention) Newton's method for the \ma equation can be written
\begin{equation}\label{eqn:newton}
-\nabla \cdot \left(B^n\nabla \varepsilon\psi\right) 
+ \nabla \cdot \left(\nabla_{\bm{x}}\left( \frac{c}{m}(\bm{x}^n)\right)\varepsilon\psi\right) 
- \varepsilon\psi \nabla_{\bm{x}} \cdot \nabla \left(\frac{c}{m}(\bm{x}^n)\right)
- |I+H(\phi^{n})| 
+ \frac{c}{m(\bm{x}^{n})} 
=  0.
\end{equation}

Equation \eqref{eqn:newton} is what we actually compute with using the finite volume software OpenFOAM.

Note that $\nabla_x (\frac{c}{m})$ has to be computed on the physical mesh $\mathcal{T}_p$: when this mesh becomes non-orthogonal this can lead to errors if we were to compute the gradient using the normal directions to each cell as in \eqref{maine}. Instead, for this gradient calculation, we follow \citet{Weller2014} and compute the gradient in the direction $\bm{d}_i$ that goes from cell centre to cell centre. Consider cell $i$ of the physical mesh with faces indexed by $f$ and neighbour cells indexed by $N_f$.

For each face, $f$, $\bm{d}_f$ is the vector from the centre of cell $i$, $\bm{x}_i$, to the centre of cell $N_f$, $\bm{x}_{N_f}$, i.e. $\bm{d}_f = \bm{x}_{N_f}-\bm{x}_i$. The gradient, $\nabla_x (\frac{c}{m})$ is calculated for cell $i$ as:
\begin{equation}
\nabla_x \left(\frac{c}{m}\right) = \left( \sum_{f\in i} \bm{d}_f \bm{d}_f^T\right)^{-1} \sum_{f\in i} \bm{d}_f \left(\left(\frac{c}{m}\right)_{N_f} - \left(\frac{c}{m}\right)_i \right).
\end{equation}

\subsection{Solving the tensorial diffusion coefficient Poisson equation}\label{solve_matrix}
In order to find the update for both the AFP method and Newton's method we must solve an equation of the form 
\begin{equation}\label{eqn:matrix1}
\nabla \cdot B \nabla \Psi  = b
\end{equation}
which gives rise to a matrix equation
\begin{equation}
M\bm{x} = \bm{b}.
\end{equation}
To compute the entries of the matrix $M$, we first split the gradient of $\Psi$ into the components normal and tangential to cell faces such that
\begin{equation}
\nabla\Psi=\hat{\mathbf{S}}_{f}\nabla_{nf}\Psi+\left((\nabla\Psi)_{f}-\left((\nabla\Psi)_{f}\cdot\hat{\mathbf{S}}_{f}\right)\hat{\mathbf{S}}_{f}\right)\label{eq:gradPsiSplit}
\end{equation}
where $\mathbf{S}_{f}$ is the face area vector; the vector normal
to each face with magnitude equal to the face area, $\hat{\mathbf{S}}_{f}$
is the unit normal vector to each face, $\nabla_{nf}\Psi$ is the gradient
of $\Psi$ calculated in direction $\hat{\mathbf{S}}_{f}$ using \eqref{maine} and $\nabla\Psi$
is the gradient of $\Psi$ calculated at cell centres using the
divergence theorem (similarly to \eqref{eqn:div_thm}) and the $()_{f}$ notation means that it is then
interpolated from cell centres to faces.

We treat implicitly the normal components of the gradient and explicitly the tangential components of the gradient.
Substituting \eqref{eq:gradPsiSplit} into \eqref{eqn:matrix1}, using the linearity of the divergence operator
and re-arranging so that implicitly calculated parts are on the left hand side of the equation gives the linear equation:
\begin{equation}
\nabla \cdot \left(B \hat{\mathbf{S}}_{f}\nabla_{n}\Psi\right) = b -  \nabla \cdot \left(B \left((\nabla\Psi)_{f}-\left((\nabla\Psi)_{f}\cdot\hat{\mathbf{S}}_{f}\right)\hat{\mathbf{S}}_{f}\right)\right).
\end{equation}
and solve this iteratively, updating the right hand side based on the current estimate of the solution. 
We have found that 3 iterations of this is sufficient to make the outer iteration process of the AFP and Newton method stable and robust. With fewer iterations the number of fixed point (or Newton) iterations to solve the \ma equation increased, and with more the number of fixed point iterations did not notably change.

\section{Numerical experiments}\label{sec:num}
\subsection{Test cases}
Meshes are generated to equidistribute a static monitor function in two dimensions on the domain $\left[-\tfrac{1}{2},\tfrac{1}{2}\right]^2$. We follow \citet{Budd2015} and \citet{Weller2016} in using a radially symmetric monitor function of the form
\begin{equation}
\label{monitor_eqn}
m(\bm{x})= 1+\alpha_1\textrm{sech}^2\left(\alpha_2\left(|\bm{x}|^2-\alpha_3^2\right)\right).
\end{equation}
Using this form we define two different monitor functions: the \emph{ring} function, where $\alpha_1= 10$, $\alpha_2 = 200$ and $\alpha_3=0.25$; and the \emph{bell} function, where $\alpha_1= 50$, $\alpha_2 = 100$ and $\alpha_3=0$. We solve the \ma equation with periodic boundary conditions in both the $x$ and $y$ directions.

For each test problem, the computational mesh $\mathcal{T}_c$ is a uniformly spaced quadrilateral grid as depicted in Figure \ref{fig:example_comp}. The number of cells in each example will vary and be denoted $N\times N$.

\subsection{Numerical linear algebra details}\label{sec:exp}
It is well known that the Laplace operator on a periodic domain has a zero eigenvalue corresponding to the constant eigenfunction. Therefore we must remove the kernel of both $\nabla^2$ in the FP method and $\nabla \cdot \left(B^n \nabla\right)$ in the AFP method. We do so by setting a fixed reference value for one of the points in the solution -- the choice of which point and what value do not influence the results.

For the Newton method, we have an added advection term in the implicit equation to be solved at each timestep. At each iteration it is an advection-diffusion equation with the \ma equation as a source term. Thus each iteration itself is hyperbolic. The boundary conditions of the hyperbolic part of the equation need to be set along every characteristic. This is not feasible as the characteristics change with every iteration.
To have a well-posed numerical problem at each iteration we impose that the solution integrates to zero globally. This is done by solving the equation
\begin{equation}\label{eqn:newtondelta}
\delta\varepsilon\psi + \nabla \cdot \left(B^n\nabla \varepsilon\psi\right) =  \nabla_{\bm{x}}\left( \frac{c}{m}(\bm{x}^n)\right) \cdot \nabla (\varepsilon\psi) - |I+H(\phi^{n})| + \frac{c}{m(\bm{x}^{n})},
\end{equation}
where $\delta$ is a scalar with dimensions $\text{area}^{-1}$. In the 2D examples considered in this \paper we set $\delta = \frac{t10^{-4}}{\min_i V_i}$, where $V_i$ is the volume of cell $i$ in the computational mesh and $t=$1m is the thickness of the domain in the direction normal to the solution domain. An appropriate value to use in 3D has yet to be determined.

The linear solver used is geometric-algebraic multigrid (GAMG) with a symmetric Gauss-Seidel smoother. The solver tolerance was set to $10^{-12}$ with a relative tolerance of $10^{-8}$. Details of this solver can be found in the documentation of OpenFOAM \citep{openfoam}. 

To compute the eigenvalues that appear in \eqref{eqn:eigenvalue} we use LAPACK \citep{lapack} -- we found that the eigenvalue calculation internal to OpenFOAM which attempts to find roots of the characteristic polynomial of the matrix is not robust to round-off errors.

\subsection{Diagnostics}
\begin{figure}[h]
\centering\begin{tikzpicture}
\begin{loglogaxis}[width=10cm,height=8cm,
ytick={2,3,4,5,6,7,8,9,10,20,30,40},
yticklabels={,3,,,,,,,10,20,30,40},
legend pos=north west,legend style={font=\tiny},
legend columns=2,xlabel=Iterations,ylabel=Wall-clock time (seconds),
enlargelimits=true,
scatter/classes={AFP_bell={mark=text,text mark=$\bell$,black},NEWTON_bell={mark=text,text mark=$\bell$,Dark2-8-1},PMA_bell={mark=text,text mark=$\bell$,Dark2-8-3},FP_bell={mark=text,text mark=$\bell$,Dark2-8-2},
AFP_ring={mark=o,black},NEWTON_ring={mark=o,Dark2-8-1},PMA_ring={mark=o,Dark2-8-3},FP_ring={mark=o,Dark2-8-2}}]
\addplot[scatter,only marks,mark options={line width=2pt}, scatter src=explicit symbolic] table[meta=label] {_home_pbrowne_OpenFOAM_results_iters_time.table};
\legend{AFP Bell,Newton Bell,PMA Bell,FP Bell,
AFP Ring,Newton Ring,PMA Ring,FP Ring}
\end{loglogaxis}
\end{tikzpicture}
 \caption{Iterations vs Wall-clock time for a number of different experiments on the Ring and Bell test cases. These are for a computational mesh of $60 \times 60$ cells, with the codes all running in serial. Both Newton's method and the AFP method are more costly per iteration due to the 3 solves at each iteration that we need in OpenFOAM to solve the linear system - this is a feature of our implementation and not necessarily the algorithms themselves. Newton's method is more costly per iteration than all the other methods due to the extra computations needed to calculate the advection term and is an inherent feature of that algorithm.}\label{fig:iter_vs_cpu}
\end{figure}
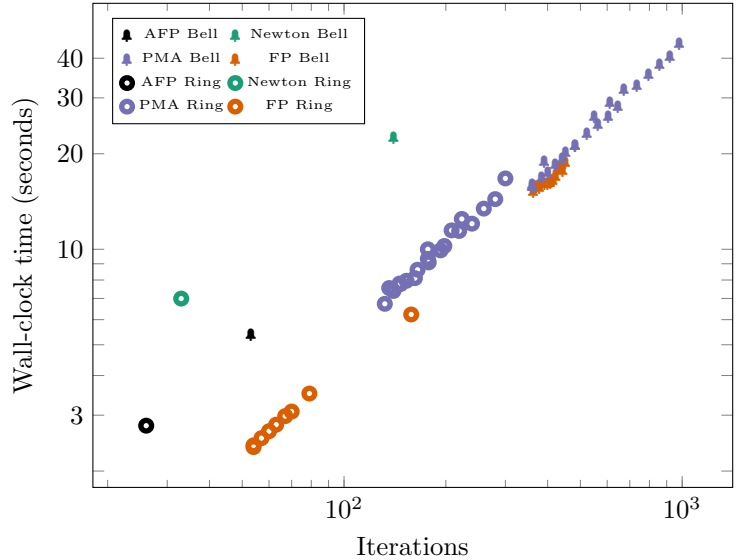

We wish to compare the efficiency of the various solution techniques in solving the \ma equation. CPU time is of course the overall goal for any numerical method, however this is not necessarily a robust measure as it will depend on both hardware and software implementation. Figure \ref{fig:iter_vs_cpu} shows, for each method, the CPU time plotted against iteration count.
We will consider the total number of outer iterations as a measure of the efficiency of the methods, and leave optimization of the codes for future investigations. 

To measure accuracy, we define an equidistribution measure following \citet{Browne2014a}: given that \eqref{eqn:ma2} says that $m(\bm{x})|I+H(\phi)| = \text{constant}$ at all points in the domain, we look at the coefficient of variation of this quantity over all the points in the mesh. That is,
\begin{equation}
\varepsilon = \frac{\sqrt{\text{Var}\{m(\bm{x})|I+H(\phi)|\}}}{\overline{m(\bm{x})|I+H(\phi)|}}.
\end{equation}
Clearly when $\varepsilon \to 0$, $m(\bm{x})|I+H(\phi)|$ is constant and thus the \ma equation is satisfied.

For all the test cases presented in this \paper, we iterate each method until the equidistribution $\varepsilon < 10^{-8}$. This is almost certainly an unnecessarily tight tolerance for mesh generation, however it will illustrate the convergence of each method to high accuracy. 

\subsection{Results}
\begin{figure}[h]
\centering

\begin{subfigure}{0.45\textwidth}
\includegraphics[width=\textwidth]{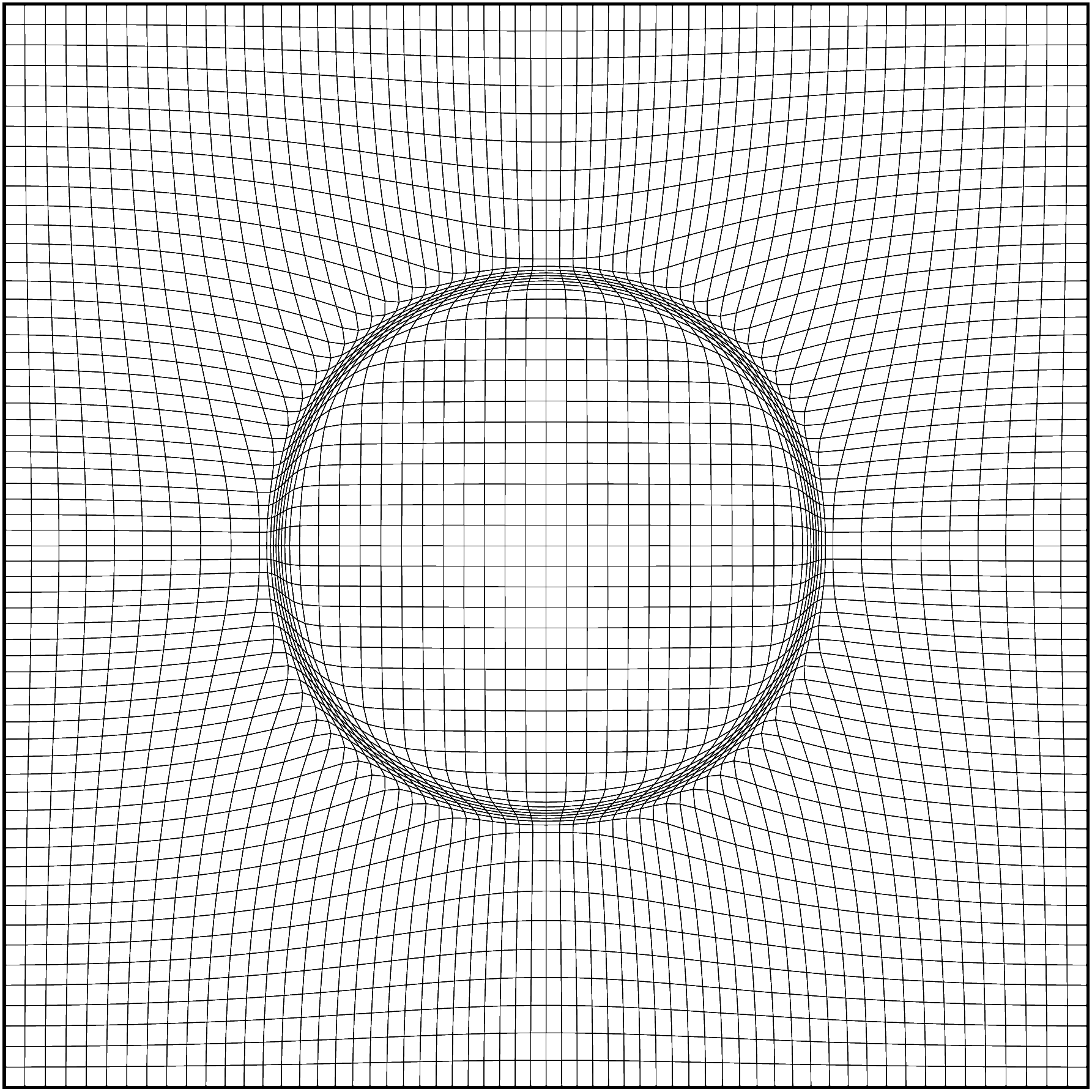}
\caption{Grid generated with AFP method for the ring test case}
\end{subfigure}\qquad
\begin{subfigure}{0.45\textwidth}
\includegraphics[width=\textwidth]{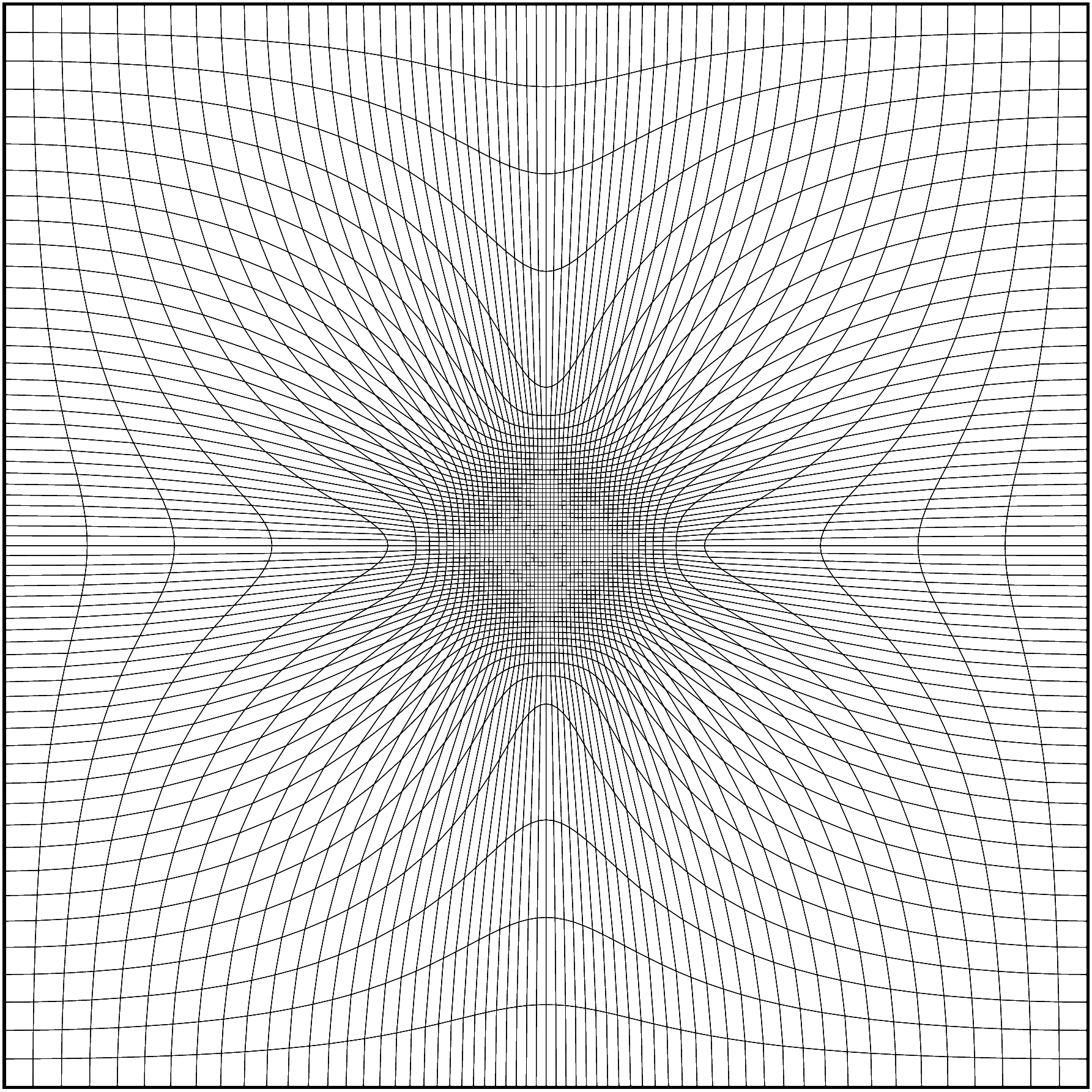}
\caption{Grid generated with AFP method for the bell test case}
\end{subfigure}
\caption{Resulting grids for the two different test cases}\label{fig:grids}
\end{figure}

Figure \ref{fig:grids} shows the resulting meshes for both the ring and the bell test case, when using the AFP method with $N=60$. At this resolution, the resulting meshes for all converged solution techniques appear identical to the eye. This is unsurprising given the unique solution to the optimal transport problem and the very tight tolerance to which we are solving the \ma equation with each method, and the fact that they all use the same spatial discretisation.

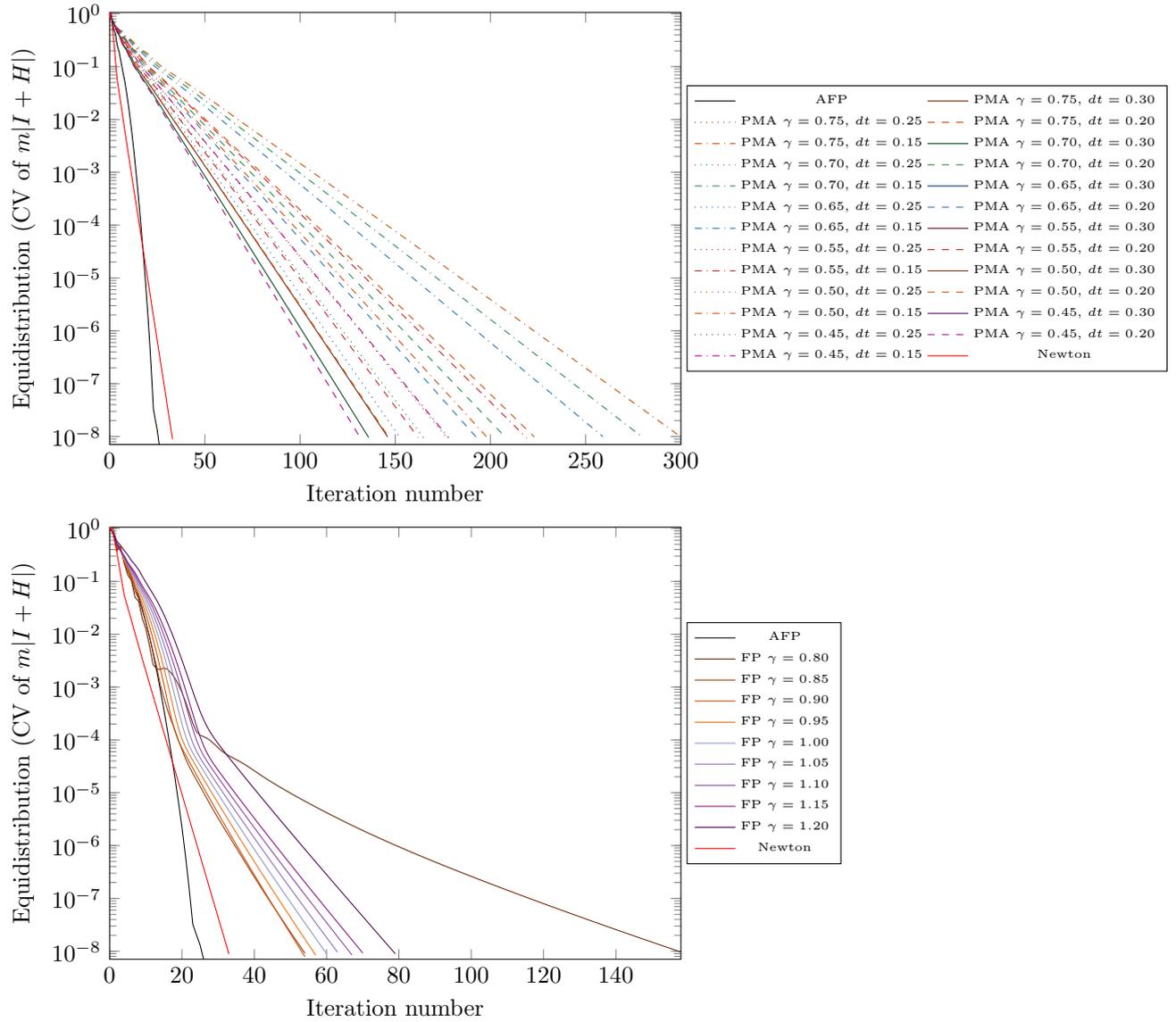
\begin{figure}
\begin{tikzpicture}
\begin{semilogyaxis}[width=10cm,height=8cm,legend style={at={(1.01,0.5)},anchor=west,font=\tiny},legend columns=2,xlabel=Iteration number,ylabel=Equidistribution (CV of $m|I+H|$),enlargelimits=false]
\addplot [black,mark=.] table[x expr=\coordindex, y index={0}] {_home_pbrowne_OpenFOAM_results_ring_AFP_60_equi};\addlegendentry{AFP};

\addplot [YlOrBr-9-9,mark=.,solid] table[x expr=\coordindex, y index={0}] {_home_pbrowne_OpenFOAM_results_ring_PMA2D_60_0.75_0.30_equi};
\addlegendentry{PMA $\gamma=0.75$, $dt = 0.30$};
\addplot [YlOrBr-9-8,mark=.,dotted] table[x expr=\coordindex, y index={0}] {_home_pbrowne_OpenFOAM_results_ring_PMA2D_60_0.75_0.25_equi};
\addlegendentry{PMA $\gamma=0.75$, $dt = 0.25$};
\addplot [YlOrBr-9-7,mark=.,dashed] table[x expr=\coordindex, y index={0}] {_home_pbrowne_OpenFOAM_results_ring_PMA2D_60_0.75_0.20_equi};
\addlegendentry{PMA $\gamma=0.75$, $dt = 0.20$};
\addplot [YlOrBr-9-7,mark=.,dashdotted] table[x expr=\coordindex, y index={0}] {_home_pbrowne_OpenFOAM_results_ring_PMA2D_60_0.75_0.15_equi};
\addlegendentry{PMA $\gamma=0.75$, $dt = 0.15$};

\addplot [Greens-9-9,mark=.,solid] table[x expr=\coordindex, y index={0}] {_home_pbrowne_OpenFOAM_results_ring_PMA2D_60_0.70_0.30_equi};
\addlegendentry{PMA $\gamma=0.70$, $dt = 0.30$};
\addplot [Greens-9-8,mark=.,dotted] table[x expr=\coordindex, y index={0}] {_home_pbrowne_OpenFOAM_results_ring_PMA2D_60_0.70_0.25_equi};
\addlegendentry{PMA $\gamma=0.70$, $dt = 0.25$};
\addplot [Greens-9-7,mark=.,dashed] table[x expr=\coordindex, y index={0}] {_home_pbrowne_OpenFOAM_results_ring_PMA2D_60_0.70_0.20_equi};
\addlegendentry{PMA $\gamma=0.70$, $dt = 0.20$};
\addplot [Greens-9-7,mark=.,dashdotted] table[x expr=\coordindex, y index={0}] {_home_pbrowne_OpenFOAM_results_ring_PMA2D_60_0.70_0.15_equi};
\addlegendentry{PMA $\gamma=0.70$, $dt = 0.15$};

\addplot [Blues-9-9,mark=.,solid] table[x expr=\coordindex, y index={0}] {_home_pbrowne_OpenFOAM_results_ring_PMA2D_60_0.65_0.30_equi};
\addlegendentry{PMA $\gamma=0.65$, $dt = 0.30$};
\addplot [Blues-9-8,mark=.,dotted] table[x expr=\coordindex, y index={0}] {_home_pbrowne_OpenFOAM_results_ring_PMA2D_60_0.65_0.25_equi};
\addlegendentry{PMA $\gamma=0.65$, $dt = 0.25$};
\addplot [Blues-9-7,mark=.,dashed] table[x expr=\coordindex, y index={0}] {_home_pbrowne_OpenFOAM_results_ring_PMA2D_60_0.65_0.20_equi};
\addlegendentry{PMA $\gamma=0.65$, $dt = 0.20$};
\addplot [Blues-9-7,mark=.,dashdotted] table[x expr=\coordindex, y index={0}] {_home_pbrowne_OpenFOAM_results_ring_PMA2D_60_0.65_0.15_equi};
\addlegendentry{PMA $\gamma=0.65$, $dt = 0.15$};

\addplot [Reds-9-9,mark=.,solid] table[x expr=\coordindex, y index={0}] {_home_pbrowne_OpenFOAM_results_ring_PMA2D_60_0.55_0.30_equi};
\addlegendentry{PMA $\gamma=0.55$, $dt = 0.30$};
\addplot [Reds-9-8,mark=.,dotted] table[x expr=\coordindex, y index={0}] {_home_pbrowne_OpenFOAM_results_ring_PMA2D_60_0.55_0.25_equi};
\addlegendentry{PMA $\gamma=0.55$, $dt = 0.25$};
\addplot [Reds-9-7,mark=.,dashed] table[x expr=\coordindex, y index={0}] {_home_pbrowne_OpenFOAM_results_ring_PMA2D_60_0.55_0.20_equi};
\addlegendentry{PMA $\gamma=0.55$, $dt = 0.20$};
\addplot [Reds-9-7,mark=.,dashdotted] table[x expr=\coordindex, y index={0}] {_home_pbrowne_OpenFOAM_results_ring_PMA2D_60_0.55_0.15_equi};
\addlegendentry{PMA $\gamma=0.55$, $dt = 0.15$};

\addplot [Oranges-9-9,mark=.,solid] table[x expr=\coordindex, y index={0}] {_home_pbrowne_OpenFOAM_results_ring_PMA2D_60_0.50_0.30_equi};
\addlegendentry{PMA $\gamma=0.50$, $dt = 0.30$};
\addplot [Oranges-9-8,mark=.,dotted] table[x expr=\coordindex, y index={0}] {_home_pbrowne_OpenFOAM_results_ring_PMA2D_60_0.50_0.25_equi};
\addlegendentry{PMA $\gamma=0.50$, $dt = 0.25$};
\addplot [Oranges-9-7,mark=.,dashed] table[x expr=\coordindex, y index={0}] {_home_pbrowne_OpenFOAM_results_ring_PMA2D_60_0.50_0.20_equi};
\addlegendentry{PMA $\gamma=0.50$, $dt = 0.20$};
\addplot [Oranges-9-7,mark=.,dashdotted] table[x expr=\coordindex, y index={0}] {_home_pbrowne_OpenFOAM_results_ring_PMA2D_60_0.50_0.15_equi};
\addlegendentry{PMA $\gamma=0.50$, $dt = 0.15$};

\addplot [RdPu-9-9,mark=.,solid] table[x expr=\coordindex, y index={0}] {_home_pbrowne_OpenFOAM_results_ring_PMA2D_60_0.45_0.30_equi};
\addlegendentry{PMA $\gamma=0.45$, $dt = 0.30$};
\addplot [RdPu-9-8,mark=.,dotted] table[x expr=\coordindex, y index={0}] {_home_pbrowne_OpenFOAM_results_ring_PMA2D_60_0.45_0.25_equi};
\addlegendentry{PMA $\gamma=0.45$, $dt = 0.25$};
\addplot [RdPu-9-7,mark=.,dashed] table[x expr=\coordindex, y index={0}] {_home_pbrowne_OpenFOAM_results_ring_PMA2D_60_0.45_0.20_equi};
\addlegendentry{PMA $\gamma=0.45$, $dt = 0.20$};
\addplot [RdPu-9-7,mark=.,dashdotted] table[x expr=\coordindex, y index={0}] {_home_pbrowne_OpenFOAM_results_ring_PMA2D_60_0.45_0.15_equi};
\addlegendentry{PMA $\gamma=0.45$, $dt = 0.15$};

\addplot [red,mark=.,solid] table[x expr=\coordindex, y index={0}] {_home_pbrowne_OpenFOAM_results_ring_Newton2d-vectorGradc_m_reconstruct_60_equi};\addlegendentry{Newton};
\end{semilogyaxis}
\end{tikzpicture}

\begin{tikzpicture}
\begin{semilogyaxis}[width=10cm,height=8cm,legend style={at={(1.01,0.5)},anchor=west,font=\tiny},xlabel=Iteration number,ylabel=Equidistribution (CV of $m|I+H|$),enlargelimits=false]
\addplot [black,mark=.] table[x expr=\coordindex, y index={0}] {_home_pbrowne_OpenFOAM_results_ring_AFP_60_equi};\addlegendentry{AFP};
\addplot [YlOrBr-9-9,mark=.] table[x expr=\coordindex, y index={0}] {_home_pbrowne_OpenFOAM_results_ring_FP2D_60_0.80_equi};\addlegendentry{FP $\gamma=0.80$};
\addplot [YlOrBr-9-8,mark=.] table[x expr=\coordindex, y index={0}] {_home_pbrowne_OpenFOAM_results_ring_FP2D_60_0.85_equi};\addlegendentry{FP $\gamma=0.85$};
\addplot [YlOrBr-9-7,mark=.] table[x expr=\coordindex, y index={0}] {_home_pbrowne_OpenFOAM_results_ring_FP2D_60_0.90_equi};\addlegendentry{FP $\gamma=0.90$};
\addplot [YlOrBr-9-6,mark=.] table[x expr=\coordindex, y index={0}] {_home_pbrowne_OpenFOAM_results_ring_FP2D_60_0.95_equi};\addlegendentry{FP $\gamma=0.95$};
\addplot [BuPu-9-5,mark=.] table[x expr=\coordindex, y index={0}] {_home_pbrowne_OpenFOAM_results_ring_FP2D_60_1.00_equi};\addlegendentry{FP $\gamma=1.00$};
\addplot [BuPu-9-6,mark=.] table[x expr=\coordindex, y index={0}] {_home_pbrowne_OpenFOAM_results_ring_FP2D_60_1.05_equi};\addlegendentry{FP $\gamma=1.05$};
\addplot [BuPu-9-7,mark=.] table[x expr=\coordindex, y index={0}] {_home_pbrowne_OpenFOAM_results_ring_FP2D_60_1.10_equi};\addlegendentry{FP $\gamma=1.10$};
\addplot [BuPu-9-8,mark=.] table[x expr=\coordindex, y index={0}] {_home_pbrowne_OpenFOAM_results_ring_FP2D_60_1.15_equi};\addlegendentry{FP $\gamma=1.15$};
\addplot [BuPu-9-9,mark=.] table[x expr=\coordindex, y index={0}] {_home_pbrowne_OpenFOAM_results_ring_FP2D_60_1.20_equi};\addlegendentry{FP $\gamma=1.20$};

\addplot [red,mark=.,solid] table[x expr=\coordindex, y index={0}] {_home_pbrowne_OpenFOAM_results_ring_Newton2d-vectorGradc_m_reconstruct_60_equi};\addlegendentry{Newton};
\end{semilogyaxis}
\end{tikzpicture}

\caption{Plots of equidistribution against iteration number for the different methods when applied to the ring test case with an initial $60\times 60$ mesh}\label{fig:equi_vs_iter_ring}
\end{figure}

\begin{figure}[h]
\begin{tikzpicture}
\begin{semilogyaxis}[width=10cm,height=8cm,legend style={at={(1.01,0.5)},anchor=west,font=\tiny},legend columns=2,xlabel=Iteration number,ylabel=Equidistribution (CV of $m|I+H|$),enlargelimits=false]
\addplot [black,mark=.] table[x expr=\coordindex, y index={0}] {_home_pbrowne_OpenFOAM_results_bell_AFP_60_equi};\addlegendentry{AFP};

\addplot [YlOrBr-9-9,mark=.,solid] table[x expr=\coordindex, y index={0}] {_home_pbrowne_OpenFOAM_results_bell_PMA2D_60_0.75_0.25_equi};
\addlegendentry{PMA $\gamma=0.75$, $dt = 0.25$};
\addplot [YlOrBr-9-8,mark=.,dotted] table[x expr=\coordindex, y index={0}] {_home_pbrowne_OpenFOAM_results_bell_PMA2D_60_0.75_0.20_equi};
\addlegendentry{PMA $\gamma=0.75$, $dt = 0.20$};
\addplot [YlOrBr-9-7,mark=.,dashed] table[x expr=\coordindex, y index={0}] {_home_pbrowne_OpenFOAM_results_bell_PMA2D_60_0.75_0.15_equi};
\addlegendentry{PMA $\gamma=0.75$, $dt = 0.15$};
\addplot [YlOrBr-9-7,mark=.,dashdotted] table[x expr=\coordindex, y index={0}] {_home_pbrowne_OpenFOAM_results_bell_PMA2D_60_0.75_0.10_equi};
\addlegendentry{PMA $\gamma=0.75$, $dt = 0.10$};

\addplot [Greens-9-9,mark=.,solid] table[x expr=\coordindex, y index={0}] {_home_pbrowne_OpenFOAM_results_bell_PMA2D_60_0.70_0.25_equi};
\addlegendentry{PMA $\gamma=0.70$, $dt = 0.25$};
\addplot [Greens-9-8,mark=.,dotted] table[x expr=\coordindex, y index={0}] {_home_pbrowne_OpenFOAM_results_bell_PMA2D_60_0.70_0.20_equi};
\addlegendentry{PMA $\gamma=0.70$, $dt = 0.20$};
\addplot [Greens-9-7,mark=.,dashed] table[x expr=\coordindex, y index={0}] {_home_pbrowne_OpenFOAM_results_bell_PMA2D_60_0.70_0.15_equi};
\addlegendentry{PMA $\gamma=0.70$, $dt = 0.15$};
\addplot [Greens-9-7,mark=.,dashdotted] table[x expr=\coordindex, y index={0}] {_home_pbrowne_OpenFOAM_results_bell_PMA2D_60_0.70_0.10_equi};
\addlegendentry{PMA $\gamma=0.70$, $dt = 0.10$};

\addplot [Blues-9-9,mark=.,solid] table[x expr=\coordindex, y index={0}] {_home_pbrowne_OpenFOAM_results_bell_PMA2D_60_0.65_0.25_equi};
\addlegendentry{PMA $\gamma=0.65$, $dt = 0.25$};
\addplot [Blues-9-8,mark=.,dotted] table[x expr=\coordindex, y index={0}] {_home_pbrowne_OpenFOAM_results_bell_PMA2D_60_0.65_0.20_equi};
\addlegendentry{PMA $\gamma=0.65$, $dt = 0.20$};
\addplot [Blues-9-7,mark=.,dashed] table[x expr=\coordindex, y index={0}] {_home_pbrowne_OpenFOAM_results_bell_PMA2D_60_0.65_0.15_equi};
\addlegendentry{PMA $\gamma=0.65$, $dt = 0.15$};
\addplot [Blues-9-7,mark=.,dashdotted] table[x expr=\coordindex, y index={0}] {_home_pbrowne_OpenFOAM_results_bell_PMA2D_60_0.65_0.10_equi};
\addlegendentry{PMA $\gamma=0.65$, $dt = 0.10$};

\addplot [Reds-9-9,mark=.,solid] table[x expr=\coordindex, y index={0}] {_home_pbrowne_OpenFOAM_results_bell_PMA2D_60_0.55_0.25_equi};
\addlegendentry{PMA $\gamma=0.55$, $dt = 0.25$};
\addplot [Reds-9-8,mark=.,dotted] table[x expr=\coordindex, y index={0}] {_home_pbrowne_OpenFOAM_results_bell_PMA2D_60_0.55_0.20_equi};
\addlegendentry{PMA $\gamma=0.55$, $dt = 0.20$};
\addplot [Reds-9-7,mark=.,dashed] table[x expr=\coordindex, y index={0}] {_home_pbrowne_OpenFOAM_results_bell_PMA2D_60_0.55_0.15_equi};
\addlegendentry{PMA $\gamma=0.55$, $dt = 0.15$};
\addplot [Reds-9-7,mark=.,dashdotted] table[x expr=\coordindex, y index={0}] {_home_pbrowne_OpenFOAM_results_bell_PMA2D_60_0.55_0.10_equi};
\addlegendentry{PMA $\gamma=0.55$, $dt = 0.10$};

\addplot [Oranges-9-9,mark=.,solid] table[x expr=\coordindex, y index={0}] {_home_pbrowne_OpenFOAM_results_bell_PMA2D_60_0.50_0.25_equi};
\addlegendentry{PMA $\gamma=0.50$, $dt = 0.25$};
\addplot [Oranges-9-8,mark=.,dotted] table[x expr=\coordindex, y index={0}] {_home_pbrowne_OpenFOAM_results_bell_PMA2D_60_0.50_0.20_equi};
\addlegendentry{PMA $\gamma=0.50$, $dt = 0.20$};
\addplot [Oranges-9-7,mark=.,dashed] table[x expr=\coordindex, y index={0}] {_home_pbrowne_OpenFOAM_results_bell_PMA2D_60_0.50_0.15_equi};
\addlegendentry{PMA $\gamma=0.50$, $dt = 0.15$};
\addplot [Oranges-9-7,mark=.,dashdotted] table[x expr=\coordindex, y index={0}] {_home_pbrowne_OpenFOAM_results_bell_PMA2D_60_0.50_0.10_equi};
\addlegendentry{PMA $\gamma=0.50$, $dt = 0.10$};

\addplot [RdPu-9-9,mark=.,solid] table[x expr=\coordindex, y index={0}] {_home_pbrowne_OpenFOAM_results_bell_PMA2D_60_0.45_0.25_equi};
\addlegendentry{PMA $\gamma=0.45$, $dt = 0.25$};
\addplot [RdPu-9-8,mark=.,dotted] table[x expr=\coordindex, y index={0}] {_home_pbrowne_OpenFOAM_results_bell_PMA2D_60_0.45_0.20_equi};
\addlegendentry{PMA $\gamma=0.45$, $dt = 0.20$};
\addplot [RdPu-9-7,mark=.,dashed] table[x expr=\coordindex, y index={0}] {_home_pbrowne_OpenFOAM_results_bell_PMA2D_60_0.45_0.15_equi};
\addlegendentry{PMA $\gamma=0.45$, $dt = 0.15$};
\addplot [RdPu-9-7,mark=.,dashdotted] table[x expr=\coordindex, y index={0}] {_home_pbrowne_OpenFOAM_results_bell_PMA2D_60_0.45_0.10_equi};
\addlegendentry{PMA $\gamma=0.45$, $dt = 0.10$};

\addplot [red,mark=.,solid] table[x expr=\coordindex, y index={0}] {_home_pbrowne_OpenFOAM_results_bell_Newton2d-vectorGradc_m_reconstruct_60_equi};\addlegendentry{Newton};
\end{semilogyaxis}
\end{tikzpicture}

\begin{tikzpicture}
\begin{semilogyaxis}[width=10cm,height=8cm,legend style={at={(1.01,0.5)},anchor=west,font=\tiny},xlabel=Iteration number,ylabel=Equidistribution (CV of $m|I+H|$),enlargelimits=false]
\addplot [black,mark=.] table[x expr=\coordindex, y index={0}] {_home_pbrowne_OpenFOAM_results_bell_AFP_60_equi};\addlegendentry{AFP};
\addplot [YlOrBr-9-9,mark=.] table[x expr=\coordindex, y index={0}] {_home_pbrowne_OpenFOAM_results_bell_FP2D_60_2.5474_equi};\addlegendentry{FP $\gamma=2.5474$};
\addplot [YlOrBr-9-8,mark=.] table[x expr=\coordindex, y index={0}] {_home_pbrowne_OpenFOAM_results_bell_FP2D_60_2.55_equi};\addlegendentry{FP $\gamma=2.55$};
\addplot [YlOrBr-9-7,mark=.] table[x expr=\coordindex, y index={0}] {_home_pbrowne_OpenFOAM_results_bell_FP2D_60_2.60_equi};\addlegendentry{FP $\gamma=2.60$};
\addplot [YlOrBr-9-6,mark=.] table[x expr=\coordindex, y index={0}] {_home_pbrowne_OpenFOAM_results_bell_FP2D_60_2.65_equi};\addlegendentry{FP $\gamma=2.65$};
\addplot [YlOrBr-9-5,mark=.] table[x expr=\coordindex, y index={0}] {_home_pbrowne_OpenFOAM_results_bell_FP2D_60_2.70_equi};\addlegendentry{FP $\gamma=2.70$};
\addplot [YlOrBr-9-4,mark=.] table[x expr=\coordindex, y index={0}] {_home_pbrowne_OpenFOAM_results_bell_FP2D_60_2.75_equi};\addlegendentry{FP $\gamma=2.75$};
\addplot [BuPu-9-3,mark=.] table[x expr=\coordindex, y index={0}] {_home_pbrowne_OpenFOAM_results_bell_FP2D_60_2.80_equi};\addlegendentry{FP $\gamma=2.80$};
\addplot [BuPu-9-4,mark=.] table[x expr=\coordindex, y index={0}] {_home_pbrowne_OpenFOAM_results_bell_FP2D_60_2.85_equi};\addlegendentry{FP $\gamma=2.85$};
\addplot [BuPu-9-5,mark=.] table[x expr=\coordindex, y index={0}] {_home_pbrowne_OpenFOAM_results_bell_FP2D_60_2.90_equi};\addlegendentry{FP $\gamma=2.90$};
\addplot [BuPu-9-6,mark=.] table[x expr=\coordindex, y index={0}] {_home_pbrowne_OpenFOAM_results_bell_FP2D_60_2.95_equi};\addlegendentry{FP $\gamma=2.95$};
\addplot [BuPu-9-7,mark=.] table[x expr=\coordindex, y index={0}] {_home_pbrowne_OpenFOAM_results_bell_FP2D_60_3.00_equi};\addlegendentry{FP $\gamma=3.00$};
\addplot [BuPu-9-8,mark=.] table[x expr=\coordindex, y index={0}] {_home_pbrowne_OpenFOAM_results_bell_FP2D_60_3.05_equi};\addlegendentry{FP $\gamma=3.05$};
\addplot [BuPu-9-9,mark=.] table[x expr=\coordindex, y index={0}] {_home_pbrowne_OpenFOAM_results_bell_FP2D_60_3.10_equi};\addlegendentry{FP $\gamma=3.10$};

\addplot [red,mark=.,solid] table[x expr=\coordindex, y index={0}] {_home_pbrowne_OpenFOAM_results_bell_Newton2d-vectorGradc_m_reconstruct_60_equi};\addlegendentry{Newton};

\end{semilogyaxis}
\end{tikzpicture}

\caption{Plots of equidistribution against iteration number for the different methods when applied to the bell test case with an initial $60\times 60$ mesh}\label{fig:equi_vs_iter_bell}
\end{figure}
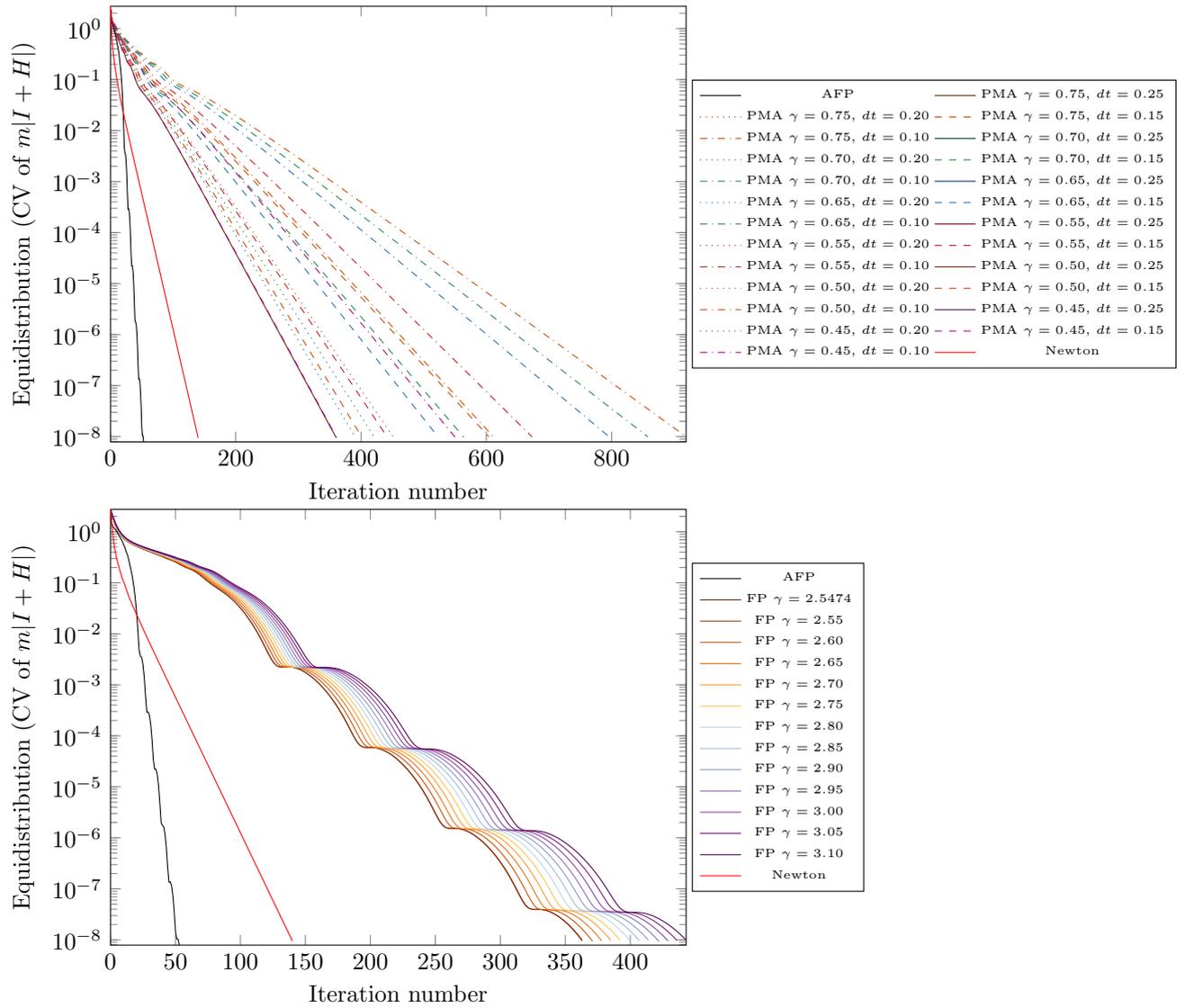

It can be seen from Figures \ref{fig:equi_vs_iter_ring} and \ref{fig:equi_vs_iter_bell} that both the fixed point iterations and the PMA technique exhibit linear convergence. However the rate of this convergence is dependent on the choice of parameters used. The AFP method, with no free parameters also shows linear convergence but at a must faster rate. Newton's method is seen to have the best performance initially when the solution is far from converged. This is due to the advection term that accelerates the solution by moving the mesh points up the gradient of the monitor function.

\clearpage
\subsection{Scaling of the methods with mesh size}

\begin{figure}[h!]
\centering
\begin{subfigure}{0.8\textwidth}
\caption{Plot of number of iterations taken against mesh size for the ring test case $\downarrow$}\label{fig:scal_ring}
\begin{tikzpicture}
\begin{axis}[width=12cm,height=8cm,legend style={at={(1.01,0.5),font=\tiny},anchor=west},xlabel=Mesh size $N$,ylabel=Iterations]
\addplot [blue,mark=o,legend entry=AFP] table [x={N}, y={iterations}] {_home_pbrowne_OpenFOAM_results_ring_AFP_data.txt};
\addplot [BuPu-9-4,mark=asterisk,legend entry={FP $\gamma=0.80$}] table [x={N}, y={iterations}] {_home_pbrowne_OpenFOAM_results_ring_FP2D_data0.80.txt};\addplot [BuPu-9-5,mark=asterisk,legend entry={FP $\gamma=0.90$}] table [x={N}, y={iterations}] {_home_pbrowne_OpenFOAM_results_ring_FP2D_data0.90.txt};\addplot [BuPu-9-6,mark=asterisk,legend entry={FP $\gamma=1.00$}] table [x={N}, y={iterations}] {_home_pbrowne_OpenFOAM_results_ring_FP2D_data1.00.txt};\addplot [BuPu-9-7,mark=asterisk,legend entry={FP $\gamma=1.10$}] table [x={N}, y={iterations}] {_home_pbrowne_OpenFOAM_results_ring_FP2D_data1.10.txt};\addplot [BuPu-9-8,mark=asterisk,legend entry={FP $\gamma=1.20$}] table [x={N}, y={iterations}] {_home_pbrowne_OpenFOAM_results_ring_FP2D_data1.20.txt};\addplot [Greens-9-9,mark=triangle,legend entry={PMA $\gamma=0.50$, $dt = 0.20$}] table [x={N}, y={iterations}] {_home_pbrowne_OpenFOAM_results_ring_PMA2D_data0.50_0.20.txt};
\addplot [Greens-9-9,mark=triangle,legend entry={PMA $\gamma=0.60$, $dt = 0.20$}] table [x={N}, y={iterations}] {_home_pbrowne_OpenFOAM_results_ring_PMA2D_data0.60_0.20.txt};
\addplot [Greens-9-9,mark=triangle,legend entry={PMA $\gamma=0.70$, $dt = 0.20$}] table [x={N}, y={iterations}] {_home_pbrowne_OpenFOAM_results_ring_PMA2D_data0.70_0.20.txt};

\addplot [Greens-9-8,mark=diamond,legend entry={PMA $\gamma=0.60$, $dt = 0.25$}] table [x={N}, y={iterations}] {_home_pbrowne_OpenFOAM_results_ring_PMA2D_data0.60_0.25.txt};
\addplot [Greens-9-8,mark=diamond,legend entry={PMA $\gamma=0.70$, $dt = 0.25$}] table [x={N}, y={iterations}] {_home_pbrowne_OpenFOAM_results_ring_PMA2D_data0.70_0.25.txt};

\addplot [Greens-9-7,mark=pentagon,legend entry={PMA $\gamma=0.70$, $dt = 0.30$}] table [x={N}, y={iterations}] {_home_pbrowne_OpenFOAM_results_ring_PMA2D_data0.70_0.30.txt};

\addplot [YlOrBr-9-9,mark=+,legend entry=Newton] table [x={N}, y={iterations}] {_home_pbrowne_OpenFOAM_results_ring_Newton2d-vectorGradc_m_reconstruct_newton_vector_data.txt};

\end{axis}
\end{tikzpicture}
\end{subfigure}

\centering\begin{subfigure}{0.8\textwidth}
\caption{Plot of number of iterations taken against mesh size for the bell test case $\downarrow$}\label{fig:scal_bell}
\begin{tikzpicture}
\begin{axis}[width=12cm,height=8cm,legend style={at={(1.01,0.5)},anchor=west,font=\tiny},xlabel=Mesh size $N$,ylabel=Iterations]
\addplot [blue,mark=o] table [x={N}, y={iterations}] {_home_pbrowne_OpenFOAM_results_bell_AFP_data.txt};\addlegendentry{AFP};
\addplot [BuPu-9-4,mark=asterisk] table [x={N}, y={iterations}] {_home_pbrowne_OpenFOAM_results_bell_FP2D_data2.60.txt};\addlegendentry{FP $\gamma=2.60$};
\addplot [BuPu-9-5,mark=asterisk] table [x={N}, y={iterations}] {_home_pbrowne_OpenFOAM_results_bell_FP2D_data2.70.txt};\addlegendentry{FP $\gamma=2.70$};
\addplot [BuPu-9-6,mark=asterisk] table [x={N}, y={iterations}] {_home_pbrowne_OpenFOAM_results_bell_FP2D_data2.80.txt};\addlegendentry{FP $\gamma=2.80$};
\addplot [BuPu-9-7,mark=asterisk] table [x={N}, y={iterations}] {_home_pbrowne_OpenFOAM_results_bell_FP2D_data2.90.txt};\addlegendentry{FP $\gamma=2.90$};
\addplot [BuPu-9-8,mark=asterisk] table [x={N}, y={iterations}] {_home_pbrowne_OpenFOAM_results_bell_FP2D_data3.00.txt};\addlegendentry{FP $\gamma=3.00$};
\addplot [BuPu-9-9,mark=asterisk] table [x={N}, y={iterations}] {_home_pbrowne_OpenFOAM_results_bell_FP2D_data3.10.txt};\addlegendentry{FP $\gamma=3.10$};
\addplot [Greens-9-9,mark=triangle] table [x={N}, y={iterations}] {_home_pbrowne_OpenFOAM_results_bell_PMA2D_data0.60_0.20.txt};
\addlegendentry{PMA $\gamma=0.60$, $dt = 0.20$};
\addplot [Greens-9-9,mark=triangle] table [x={N}, y={iterations}] {_home_pbrowne_OpenFOAM_results_bell_PMA2D_data0.70_0.20.txt};
\addlegendentry{PMA $\gamma=0.70$, $dt = 0.20$};

\addplot [YlOrBr-9-9,mark=+,legend entry=Newton] table [x={N}, y={iterations}] {_home_pbrowne_OpenFOAM_results_bell_Newton2d-vectorGradc_m_reconstruct_newton_vector_data.txt};
\addlegendentry{Newton};
\end{axis}
\end{tikzpicture}
\end{subfigure}
\caption{Scaling of the different algorithms as the resolution of the computational mesh increases.}\label{fig:scal}
\end{figure}
 
Figure \ref{fig:scal} shows how the number of iterations for each method varies as a function of the (square root of the) number of cells in the computational mesh varies. There are multiple lines for both the PMA method and the FP method due to the different free parameters in each technique. In contrast, there is only one line for the AFP method as it has no free parameters. We only show Newton's method for a single value of the constant $\delta$. Considering first the PMA method, we see that the number of iterations is (for large $N$) independent of the mesh size. This is consistent with previous studies \citep{Browne2014a}. However, note that for the ring test case the PMA method failed in 3 out of the 9 parameter sets: the lines $(\gamma=0.5,dt=0.25)$, $(\gamma=0.5,dt=0.3)$, and $(\gamma=0.6,dt=0.3)$ all do not appear in Figure \ref{fig:scal_ring}. Further, only 2 of the parameter sets converged for the bell test case. This highlights the importance of these free parameters for the tuning of the PMA method.

For the fixed point method, Figure \ref{fig:scal} shows that the optimal regularisation parameter $\gamma$ is strongly problem dependent. For the ring test case, a $\gamma$ in the interval $[0.8,1.2]$ was reasonable, whereas in the bell test case $\gamma \in [2.6,3.1]$ was more appropriate. For the ring test case, Figure \ref{fig:scal_ring} shows that the number of iterations is independent of $N$ as in the PMA method. The general pattern that a smaller $\gamma$ leads to fewer iterations with the FP method can be seen clearly in Figure \ref{fig:scal_ring}. There is, however, a limit to how small $\gamma$ can be before the algorithm fails. This is seen with $\gamma=0.8$, that only converges for the smallest $N=60$ case. For the bell test case (Figure \ref{fig:scal_bell}) it is clear that the optimal $\gamma$ is dependent on the mesh size $N$. The higher the resolution, the larger $\gamma$ is needed for the method to converge, at the expense of taking more iterations.

For the AFP method, Figure \ref{fig:scal} shows that there may be a slight increase in the number of iterations taken as the resolution increases. However, the order of magnitude remains the same and the method converges for every mesh size without having to choose any parameters.
\subsection{Problems with the Newton method}
Consider now the behaviour of the Newton solver shown in Figure \ref{fig:scal}. Notice immediately that the Newton method failed for the  bell test case when the problem size is greater than $150\times150$. The method fails due to catastrophic mesh tangling 
and we have investigated whether this could be caused
by
numerical errors in the calculation of  $\nabla_x \left(\frac{c}{m}\right)$ on the physical mesh.

For the form of the monitor function that we consider in this paper (see \eqref{monitor_eqn}) it is a simple exercise to use symbolic algebra to derive an analytic expression for $\nabla_x \left(\frac{c}{m}\right)$. 
Numerical tests were performed using the analytic rather than the numerical gradient of $\nabla_x \left(\frac{c}{m}\right)$. We found that using the smooth, analytic gradient
neither increased the robustness of the Newton method (i.e. achieved convergence when the numerical calculation of the gradient failed) nor increased the efficiency of the method (i.e. the results remained within a single iteration of the numerical gradient calculations).
Therefore the slow convergence or divergence of the Newton method 
appears not to be due to numerical errors in calculating $\nabla_x \left(\frac{c}{m}\right)$

One can also add under-relaxation to the Newton method by increasing the value of $\delta$ in \eqref{eqn:newtondelta} -- we found that if $\delta$ were increased to a level that make the method convergent then the number of iterations needed were prohibitively large.
It is therefore not clear why the Newton method fails in some cases or why it is less convergent than the fixed point method when both are near convergence.

\subsection{Use of the adaptive regularisation term}
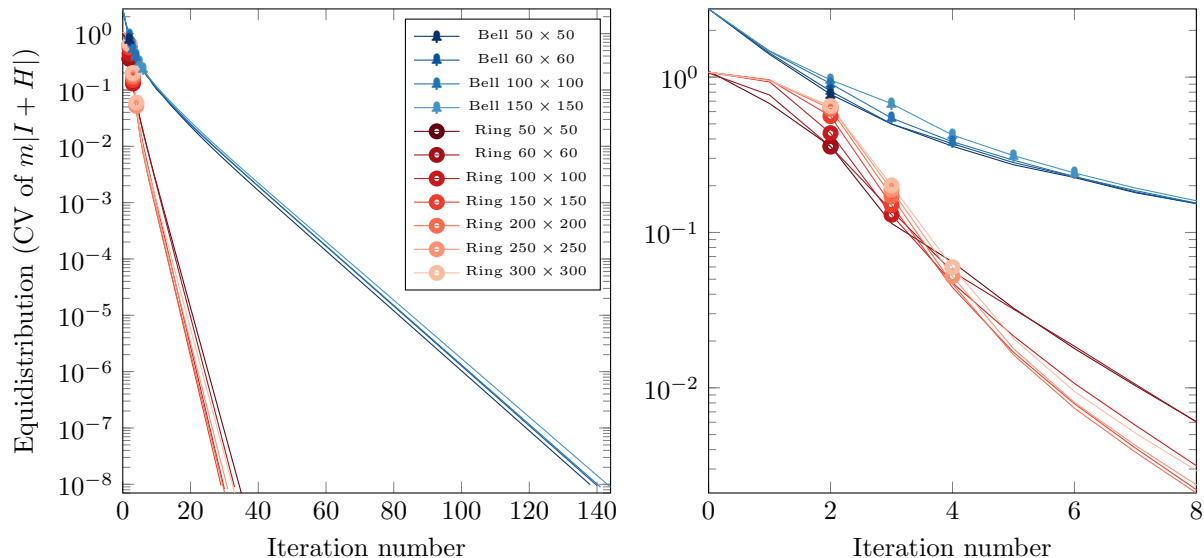
\begin{figure}[h]
\begin{tikzpicture}
\begin{semilogyaxis}[width=8cm,height=8cm,legend style={font=\tiny},legend columns=1,xlabel=Iteration number,ylabel=Equidistribution (CV of $m|I+H|$),enlargelimits=false,mark options={line width=2pt}]

\addplot [Blues-9-9,mark=.,solid,forget plot] table[x expr=\coordindex, y index={0}] {_home_pbrowne_OpenFOAM_results_bell_Newton2d-vectorGradc_m_reconstruct_50_equi};\addlegendimage{Blues-9-9,mark=text, text mark=$\bell$}\addlegendentry{Bell $50\times 50$};
\addplot [Blues-9-8,mark=.,solid,forget plot] table[x expr=\coordindex, y index={0}] {_home_pbrowne_OpenFOAM_results_bell_Newton2d-vectorGradc_m_reconstruct_60_equi};\addlegendimage{Blues-9-8,mark=text, text mark=$\bell$}\addlegendentry{Bell $60\times 60$};
\addplot [Blues-9-7,mark=.,solid,forget plot] table[x expr=\coordindex, y index={0}] {_home_pbrowne_OpenFOAM_results_bell_Newton2d-vectorGradc_m_reconstruct_100_equi};\addlegendimage{Blues-9-7,mark=text, text mark=$\bell$}\addlegendentry{Bell $100\times 100$};
\addplot [Blues-9-6,mark=.,solid,forget plot] table[x expr=\coordindex, y index={0}] {_home_pbrowne_OpenFOAM_results_bell_Newton2d-vectorGradc_m_reconstruct_150_equi};\addlegendimage{Blues-9-6,mark=text, text mark=$\bell$}\addlegendentry{Bell $150\times 150$};

\addplot [Reds-9-9,mark=.,solid,forget plot] table[x expr=\coordindex, y index={0}] {_home_pbrowne_OpenFOAM_results_ring_Newton2d-vectorGradc_m_reconstruct_50_equi};\addlegendimage{Reds-9-9,mark=o};\addlegendentry{Ring $50\times 50$};

\addplot [Reds-9-8,mark=.,solid,forget plot] table[x expr=\coordindex, y index={0}] {_home_pbrowne_OpenFOAM_results_ring_Newton2d-vectorGradc_m_reconstruct_60_equi};\addlegendimage{Reds-9-8,mark=o};\addlegendentry{Ring $60\times 60$};

\addplot [Reds-9-7,mark=.,solid,forget plot] table[x expr=\coordindex, y index={0}] {_home_pbrowne_OpenFOAM_results_ring_Newton2d-vectorGradc_m_reconstruct_100_equi};\addlegendimage{Reds-9-7,mark=o};\addlegendentry{Ring $100\times 100$};

\addplot [Reds-9-6,mark=.,solid,forget plot] table[x expr=\coordindex, y index={0}] {_home_pbrowne_OpenFOAM_results_ring_Newton2d-vectorGradc_m_reconstruct_150_equi};\addlegendimage{Reds-9-6,mark=o};\addlegendentry{Ring $150\times 150$};

\addplot [Reds-9-5,mark=.,solid,forget plot] table[x expr=\coordindex, y index={0}] {_home_pbrowne_OpenFOAM_results_ring_Newton2d-vectorGradc_m_reconstruct_200_equi};\addlegendimage{Reds-9-5,mark=o};\addlegendentry{Ring $200\times 200$};

\addplot [Reds-9-4,mark=.,solid,forget plot] table[x expr=\coordindex, y index={0}] {_home_pbrowne_OpenFOAM_results_ring_Newton2d-vectorGradc_m_reconstruct_250_equi};\addlegendimage{Reds-9-4,mark=o};\addlegendentry{Ring $250\times 250$};

\addplot [Reds-9-3,mark=.,solid,forget plot] table[x expr=\coordindex, y index={0}] {_home_pbrowne_OpenFOAM_results_ring_Newton2d-vectorGradc_m_reconstruct_300_equi};\addlegendimage{Reds-9-3,mark=o};\addlegendentry{Ring $300\times 300$};

\addplot [only marks,mark=o,Reds-9-9] table[x index={0}, y index={1}] {_home_pbrowne_OpenFOAM_results_ring_Newton2d-vectorGradc_m_reconstruct_50_regs};
\addplot [only marks,mark=o,Reds-9-8] table[x index={0}, y index={1}] {_home_pbrowne_OpenFOAM_results_ring_Newton2d-vectorGradc_m_reconstruct_60_regs};
\addplot [only marks,mark=o,Reds-9-7] table[x index={0}, y index={1}] {_home_pbrowne_OpenFOAM_results_ring_Newton2d-vectorGradc_m_reconstruct_100_regs};
\addplot [only marks,mark=o,Reds-9-6] table[x index={0}, y index={1}] {_home_pbrowne_OpenFOAM_results_ring_Newton2d-vectorGradc_m_reconstruct_150_regs};
\addplot [only marks,mark=o,Reds-9-5] table[x index={0}, y index={1}] {_home_pbrowne_OpenFOAM_results_ring_Newton2d-vectorGradc_m_reconstruct_200_regs};
\addplot [only marks,mark=o,Reds-9-4] table[x index={0}, y index={1}] {_home_pbrowne_OpenFOAM_results_ring_Newton2d-vectorGradc_m_reconstruct_250_regs};
\addplot [only marks,mark=o,Reds-9-3] table[x index={0}, y index={1}] {_home_pbrowne_OpenFOAM_results_ring_Newton2d-vectorGradc_m_reconstruct_300_regs};

\addplot [only marks,mark=text, text mark=$\bell$,Blues-9-6] table[x index={0}, y index={1}] {_home_pbrowne_OpenFOAM_results_bell_Newton2d-vectorGradc_m_reconstruct_150_regs};
\addplot [only marks,mark=text, text mark=$\bell$,Blues-9-7] table[x index={0}, y index={1}] {_home_pbrowne_OpenFOAM_results_bell_Newton2d-vectorGradc_m_reconstruct_100_regs};
\addplot [only marks,mark=text, text mark=$\bell$,Blues-9-8] table[x index={0}, y index={1}] {_home_pbrowne_OpenFOAM_results_bell_Newton2d-vectorGradc_m_reconstruct_60_regs};
\addplot [only marks,mark=text, text mark=$\bell$,Blues-9-9] table[x index={0}, y index={1}] {_home_pbrowne_OpenFOAM_results_bell_Newton2d-vectorGradc_m_reconstruct_50_regs};

\end{semilogyaxis}
\end{tikzpicture}
\begin{tikzpicture}
\begin{semilogyaxis}[width=8cm,height=8cm,xlabel=Iteration number,enlargelimits=false,xmax=8,mark options={line width=2pt}]

\addplot [Blues-9-9,mark=.,solid,forget plot] table[x expr=\coordindex, y index={0}] {_home_pbrowne_OpenFOAM_results_bell_Newton2d-vectorGradc_m_reconstruct_50_equi};\addlegendimage{Blues-9-9,mark=text, text mark=$\bell$}\addlegendentry{Bell $50\times 50$};
\addplot [Blues-9-8,mark=.,solid,forget plot] table[x expr=\coordindex, y index={0}] {_home_pbrowne_OpenFOAM_results_bell_Newton2d-vectorGradc_m_reconstruct_60_equi};\addlegendimage{Blues-9-8,mark=text, text mark=$\bell$}\addlegendentry{Bell $60\times 60$};
\addplot [Blues-9-7,mark=.,solid,forget plot] table[x expr=\coordindex, y index={0}] {_home_pbrowne_OpenFOAM_results_bell_Newton2d-vectorGradc_m_reconstruct_100_equi};\addlegendimage{Blues-9-7,mark=text, text mark=$\bell$}\addlegendentry{Bell $100\times 100$};
\addplot [Blues-9-6,mark=.,solid,forget plot] table[x expr=\coordindex, y index={0}] {_home_pbrowne_OpenFOAM_results_bell_Newton2d-vectorGradc_m_reconstruct_150_equi};\addlegendimage{Blues-9-6,mark=text, text mark=$\bell$}\addlegendentry{Bell $150\times 150$};

\addplot [Reds-9-9,mark=.,solid,forget plot] table[x expr=\coordindex, y index={0}] {_home_pbrowne_OpenFOAM_results_ring_Newton2d-vectorGradc_m_reconstruct_50_equi};\addlegendimage{Reds-9-9,mark=o};\addlegendentry{Ring $50\times 50$};

\addplot [Reds-9-8,mark=.,solid,forget plot] table[x expr=\coordindex, y index={0}] {_home_pbrowne_OpenFOAM_results_ring_Newton2d-vectorGradc_m_reconstruct_60_equi};\addlegendimage{Reds-9-8,mark=o};\addlegendentry{Ring $60\times 60$};

\addplot [Reds-9-7,mark=.,solid,forget plot] table[x expr=\coordindex, y index={0}] {_home_pbrowne_OpenFOAM_results_ring_Newton2d-vectorGradc_m_reconstruct_100_equi};\addlegendimage{Reds-9-7,mark=o};\addlegendentry{Ring $100\times 100$};

\addplot [Reds-9-6,mark=.,solid,forget plot] table[x expr=\coordindex, y index={0}] {_home_pbrowne_OpenFOAM_results_ring_Newton2d-vectorGradc_m_reconstruct_150_equi};\addlegendimage{Reds-9-6,mark=o};\addlegendentry{Ring $150\times 150$};

\addplot [Reds-9-5,mark=.,solid,forget plot] table[x expr=\coordindex, y index={0}] {_home_pbrowne_OpenFOAM_results_ring_Newton2d-vectorGradc_m_reconstruct_200_equi};\addlegendimage{Reds-9-5,mark=o};\addlegendentry{Ring $200\times 200$};

\addplot [Reds-9-4,mark=.,solid,forget plot] table[x expr=\coordindex, y index={0}] {_home_pbrowne_OpenFOAM_results_ring_Newton2d-vectorGradc_m_reconstruct_250_equi};\addlegendimage{Reds-9-4,mark=o};\addlegendentry{Ring $250\times 250$};

\addplot [Reds-9-3,mark=.,solid,forget plot] table[x expr=\coordindex, y index={0}] {_home_pbrowne_OpenFOAM_results_ring_Newton2d-vectorGradc_m_reconstruct_300_equi};\addlegendimage{Reds-9-3,mark=o};\addlegendentry{Ring $300\times 300$};

\addplot [only marks,mark=text, text mark=$\bell$,Blues-9-6] table[x index={0}, y index={1}] {_home_pbrowne_OpenFOAM_results_bell_Newton2d-vectorGradc_m_reconstruct_150_regs};
\addplot [only marks,mark=text, text mark=$\bell$,Blues-9-7] table[x index={0}, y index={1}] {_home_pbrowne_OpenFOAM_results_bell_Newton2d-vectorGradc_m_reconstruct_100_regs};
\addplot [only marks,mark=text, text mark=$\bell$,Blues-9-8] table[x index={0}, y index={1}] {_home_pbrowne_OpenFOAM_results_bell_Newton2d-vectorGradc_m_reconstruct_60_regs};
\addplot [only marks,mark=text, text mark=$\bell$,Blues-9-9] table[x index={0}, y index={1}] {_home_pbrowne_OpenFOAM_results_bell_Newton2d-vectorGradc_m_reconstruct_50_regs};

\addplot [only marks,mark=o,Reds-9-9] table[x index={0}, y index={1}] {_home_pbrowne_OpenFOAM_results_ring_Newton2d-vectorGradc_m_reconstruct_50_regs};
\addplot [only marks,mark=o,Reds-9-8] table[x index={0}, y index={1}] {_home_pbrowne_OpenFOAM_results_ring_Newton2d-vectorGradc_m_reconstruct_60_regs};
\addplot [only marks,mark=o,Reds-9-7] table[x index={0}, y index={1}] {_home_pbrowne_OpenFOAM_results_ring_Newton2d-vectorGradc_m_reconstruct_100_regs};
\addplot [only marks,mark=o,Reds-9-6] table[x index={0}, y index={1}] {_home_pbrowne_OpenFOAM_results_ring_Newton2d-vectorGradc_m_reconstruct_150_regs};
\addplot [only marks,mark=o,Reds-9-5] table[x index={0}, y index={1}] {_home_pbrowne_OpenFOAM_results_ring_Newton2d-vectorGradc_m_reconstruct_200_regs};
\addplot [only marks,mark=o,Reds-9-4] table[x index={0}, y index={1}] {_home_pbrowne_OpenFOAM_results_ring_Newton2d-vectorGradc_m_reconstruct_250_regs};
\addplot [only marks,mark=o,Reds-9-3] table[x index={0}, y index={1}] {_home_pbrowne_OpenFOAM_results_ring_Newton2d-vectorGradc_m_reconstruct_300_regs};
\legend{};\end{semilogyaxis}
\end{tikzpicture}
\caption{Plots of equidistribution against iteration number for the Newton method with the marks showing the iterations where regularisation occured (i.e. $\gamma \ne 0$ in \eqref{eqn:eigenvalue}).
On the left is show the convergence for every iteration whilst on the right we zoom in to the initial iterations where the regularisation was needed.
}\label{fig:equi_vs_iter_reg}
\end{figure}
 Figure \ref{fig:equi_vs_iter_reg} shows the iterations at which the regularisation term in equation \eqref{eqn:eigenvalue} was used. In the numerical experiments we have conducted this was only ever needed in the Newton method and not the AFP method. However, when fewer correcting iterations are used to solve the numerical linear algebra problem as discussed in Section \ref{solve_matrix}, this regularisation is also needed in the AFP method. Note from Figure \ref{fig:equi_vs_iter_reg} that $\gamma \ne 0$ only in the first few iterations, when the solution is far from converged to the unique convex solution.

\section{Conclusions}\label{sec:conc}
We have introduced a new iterative method for solving the \ma equation for mesh redistribution based on linearising the determinant of the Jacobian of the optimal transport map about a previous solution. Using this linearisation we have added a regularisation based on the spectrum of the discrete linearised operator to ensure ellipticity at each step of the iterative procedure. This has resulted in a method with no free parameters which has been robust to all the test cases we have considered.

We have shown that the method is more efficient than a parabolic relaxation of the \ma equation and is significantly more efficient than a fixed point method based on a linearisation about $0$ when such a method requires a large amount of under-relaxation to converge.

We have considered a complete linearisation of the \ma equation that leads to Newton's method to solve the \ma equation. We have shown that this method does not always lead to a reduction in the number of iterations needed to solve the \ma equation when compared to the AFP method. 
Newton's method was found to exhibit only linear convergence, and not quadratic or super-linear convergence that would make it more attractive as a solution algorithm compared with the simpler algorithms that can provide linear convergence with similar or better rates.

This \paper has considered applying different algorithms directly to given monitor functions without applying any \textit{ad hoc} smoothing (also referred to as filtering) of the monitor function. Such smoothing will necessarily lead to smoother meshes (as the regularity of the mesh is equivalent to the regularity of the monitor function), however it obscures some of the convergence properties of the solution process. The results given in this \paper may therefore be different when smoothing is present.

There is still an open question of what convergence criteria to use when stopping an algorithm for mesh redistribution. In this paper, we have used a global statistic based on the equidistribution measure defined at the centre of each finite volume cell. 
If a redistributed mesh is going to be used for 
the numerical solution of a set of PDEs,
 equidistribution may
not need to be satisfied to a tight tolerance, and the condition that it holds at a specific point in each cell is likely not necessary. Future investigations into optimally transported mesh redistribution techniques may wish to consider stopping criteria that hold only somewhere in a given cell.

Finally, note that the results stated here are valid only in Euclidean geometries (specifically on the plane). When solving mesh redistribution on the sphere or other manifolds, an analogous linearisation of the determinant of the Jacobian of the optimal transport map must take into account of the curvature of the manifold. The linearisation technique given in this \paper may provide a framework to form such a linearisation.

\section{Code availability}
The codes used to implement the test problems in this paper, using OpenFOAM, are available from \url{http://researchdata.reading.ac.uk/SOME_URL} \todo{Link available following \paper acceptance}

The codes are available directly from github: 
\url{https://github.com/AtmosFOAM/AMMM/releases/tag/v0.2}

\section{Acknowledgements}\label{sec:ack}
This work was supported by NERC grant NE/M013693/1 (PAB \& HW), EPSRC grant EP/L016613/1 (JP) and EPSRC grant EP/P000835/1 (TP).

\bibliographystyle{apalike}

\newpage
\listoftodos
\end{document}